\documentclass[journal,twoside,web]{ieeecolor}
\usepackage{generic}
\usepackage{textcomp}
\usepackage{color}

\usepackage{amsthm}
\usepackage{amsmath}
\usepackage{hyperref}
\hypersetup{hidelinks}
\usepackage[english]{babel}
\usepackage{comment}
\theoremstyle{definition}
\newtheorem{proposition}{Proposition}
\newtheorem{definition}{Definition}
\newtheorem{corollary}{Corollary}

\newtheorem{theorem}{Theorem}
\newtheorem{remark}{Remark}

\renewenvironment{proof}{{\itshape Proof.}}{\qedsymbol}

\DeclareMathOperator*{\argmax}{arg\,max}
\date{October 2024}

\newcommand{\shnote}[1]%
    {\textcolor{magenta}{ #1}}

\newcommand{\bynote}[1]%
    { \textbf{\textcolor{red}{#1}} }

\newcommand{\zgnote}[1]%
    { \textbf{\textcolor{blue}{#1}} }

\newcommand{\zgrvs}[1]%
    { \textcolor{red}{#1}}

\newcommand{\target}{\mathcal{T}}
\newcommand{\obs}{\mathcal{C}}
\newcommand{\tinit}{t}
\newcommand{\tvar}{s}
\newcommand{\dstb}{d} % Disturbance
\newcommand{\csig}{\mathrm{u}}
\newcommand{\dsig}{\mathrm{d}}

\newcommand{\dset}{\mathcal{D}}
\newcommand{\dfset}{\mathbb{D}}
\newcommand{\dmap}{\lambda} % dstb strategy
\newcommand{\Dmap}{\Lambda} % dstb strategy set

\newcommand{\R}{\mathbb{R}}
\newcommand{\cset}{\mathcal{U}} %tracker control set
\newcommand{\cfset}{\mathbb{U}} %tracker control funciton set
\newcommand{\mrcis}{\mathcal{I}_m}

\newcommand{\constraint}{\mathcal{C}}

\newcommand{\VSA}{V_{\gamma}^{\text{SA}}}
\newcommand{\clvf}{V_\theta^{\text{CLVF}}}
\newcommand{\sclvf}{\bar{V}_\theta^{\text{CLVF}}}
\newcommand{\minval}{H^\infty_m} % minimal Value 
\usepackage{booktabs}
\usepackage{graphicx}
\usepackage{booktabs}
\usepackage{amsmath} % assumes amsmath package installed
\let\labelindent\relax
\usepackage{enumitem}

\usepackage{amsfonts}
\usepackage{amssymb}  % assumes amsmath package installed
\usepackage{mathtools} %% for \coloneqq
\usepackage{mathrsfs}
\usepackage{balance}
\usepackage{url}
\usepackage{hyperref}
\usepackage{algorithm}
\usepackage{algpseudocode}
\DeclareMathAlphabet{\mathpzc}{OT1}{pzc}{m}{it}
\usepackage{makecell}
\usepackage[noadjust]{cite}
\begin{document}

\title{Solving Reach- and Stabilize-Avoid Problems Using Discounted Reachability}

\author{Boyang Li*, \IEEEmembership{Student Member, IEEE}, Zheng Gong*, \IEEEmembership{Student Member, IEEE}, and Sylvia Herbert, \IEEEmembership{Member, IEEE}
\thanks{This work is supported by ONR YIP (\#N00014-22-1-2292) and the UCSD JSOE Early Career Faculty Award.}
\thanks{*Both authors contributed equally to this work. All authors are in Mechanical and Aerospace Engineering at UC San Diego (e-mail: 
\{\href{mailto:bol025@ucsd.edu}{bol025}, {\href{mailto:zhgong@ucsd.edu}{zhgong}, \href{mailto:sherbert@ucsd.edu}{sherbert}\}@ucsd.edu.})
}
\thanks{This article has been accepted for publication in the IEEE Transactions on Automatic Control. The version of record is available at \href{https://doi.org/10.1109/TAC.2026.3693989}{IEEE Xplore} (DOI: 10.1109/TAC.2026.3693989).}
\thanks{\copyright~2026 IEEE. Personal use of this material is permitted. Permission from IEEE must be obtained for all other uses, in any current or future media, including reprinting/republishing this material for advertising or promotional purposes, creating new collective works, for resale or redistribution to servers or lists, or reuse of any copyrighted component of this work in other works.}}
\maketitle

\begin{abstract}
In this article, we consider the \textit{infinite-horizon reach‐avoid (RA) and stabilize‐avoid (SA) zero-sum game} problems for general nonlinear continuous-time systems, where the goal is to find the set of states that can be controlled to reach or stabilize to a target set, without violating constraints even under the worst-case disturbance. Based on the Hamilton-Jacobi reachability method, we address the RA problem by designing a new Lipschitz continuous RA value function, whose zero sublevel set exactly characterizes the RA set. We establish that the associated Bellman backup operator is contractive and that the RA value function is the unique viscosity solution of a Hamilton–Jacobi variational inequality. Finally, we develop a two-step framework for the SA problem by integrating our RA strategies with a recently proposed Robust Control Lyapunov-Value Function, thereby ensuring both target reachability and long-term stability. We numerically verify our RA and SA frameworks on a 3D Dubins car system to demonstrate the efficacy of the proposed approach. Code is available at \url{https://github.com/UCSD-SASLab/Infinite-Horizon-Reach-and-Stabilize-Avoid}.
\end{abstract}

%%%%%%%%%%%%%%%%%%%%%%%%%%%%%%%%%%%%%%%%%%%%%%%%%%%%%%%%%%%%%%%%%%%%%%%%%%%%%%%%%%%%%%%%%%%%%%%%%%%%%%%%%%%%%%%%%%%%%%%%%%%%%%%%%%%%%%%%%%%%%%%%%%%%%%%%%%%%

\section{INTRODUCTION}
Reach-Avoid (RA) problems have gained significant attention in recent years due to their broad range of applications in engineering, especially in controlling systems where safety or strategic decision-making is critical. It is characterized by a two-player zero-sum game, where the control as one of the players aims to steer the system into a target set while avoiding an unsafe set, and the disturbance as another player strives to prevent the first from succeeding.
%In many practical scenarios, for a given dynamic system,  it is necessary to find the set of initial states and design associated control strategies that not only ensures the system reaches a desired state but also respects a set of constraints, all while possibly facing unknown disturbances or adversarial interference. 
%With the increasingly frequent deployment of autonomous system into our daily lives, how to ensure safey operation of complex dynamic systems while satisfying some performance objectives have become a popular topic.

The finite-horizon RA game aims to identify the set of initial states (called the RA set) and the corresponding control law that can achieve such goals, i.e., reaching and avoiding, in a predefined finite time horizon.
%that can reach a target set within a prescribed finite time while avoiding an unsafe set, even subject to the worst-case disturbance within some pre-specified bound. 
For continuous-time systems, Hamilton-Jacobi (HJ) reachability analysis is a widely used approach to solve such problems \cite{Zidani_2010, Altarovici_Bokanowski_Zidani_2013, Fisac_Chen_Tomlin_Sastry_2015, Barron2018}. This method allows for the design of a value function whose gradients can be used for control synthesis, and the sign of the value function at a given state characterizes the safety and performance of that state. The computation of this value function involves solving a time-dependent HJ partial differential equation \cite{Mitchell_TV}, which typically relies on the value iteration that recursively applies a Bellman backup operator until convergence. 
%Additionally, this approach yields a closed-loop control policy in continuous-time scenarios for zero-sum games with a finite time horizons, ensuring that an agent reaches a specified target set within a finite time frame while ensuring safety.  

%\shnote{cite methods to solve reach-avoid  problems? Anything beyond standard dynamic programming (e.g. Ram/Shreyas' work)?}

%Infinite-horizon reach-avoid (RA) problem seeks to find the set of initial states and the associated controller to drive the system to a target set in finite time while avoiding undesirable states, even subject to the worst-case disturbances within pre-specified bounds. 
However, finite-horizon formulations have a major drawback. Since the value function inherently depends on time, the associated Bellman backup operator changes with the remaining time-to-go and is not a contraction mapping. As a result, the convergence of value iteration becomes sensitive to initialization, and correct initialization is required to ensure convergence to the correct solution. In contrast, a contraction mapping guarantees convergence from any arbitrary initialization and can even accelerate the computation when a good initialization is used, which are the theoretical pivots for methods like warm-starting \cite{warm-starting} and reinforcement learning \cite{Bridging} that aim for efficient computation and scalability.
%the challenge of determining an adequate time horizon for ensuring the feasibility of the zero-sum game, particularly in complex or uncertain dynamical systems.

%This limitation motivates the development of the infinite-horizon reach-avoid zero-sum game, which cares about if the game will be feasible in some finite time and thus eliminates the requirement of prescribing a time for reaching the target set. The infinite-horizon RA has only been explored in a few prior works, and mainly for discrete-time systems, where it was tackled by leveraging reinforcement learning techniques, \cite{Hsu*_Rubies-Royo*_Tomlin_Fisac_2021, hsunguyen2023isaacs, wang2024magicsadversarialrlminimax,li2024certifiable}. However, they either fail to provide deterministic safety guarantees or yield conservative RA set in practice. there has been interest in the infinite-horizon RA game 

This limitation motivates the development of the infinite-horizon RA zero-sum game.
For discrete-time systems, prior work leverages reinforcement learning to approximate RA value functions, where the convergence is guaranteed by designing a contractive Bellman backup operator that is induced from a discount factor \cite{Hsu*_Rubies-Royo*_Tomlin_Fisac_2021, HSU2023103811, li2023learningpredictivesafetyfilter, hsunguyen2023isaacs, wang2024magicsadversarialrlminimax,nguyen2024gameplay}. However, these works fail to provide deterministic safety guarantees or yield conservative RA sets in practice, although point-wise guarantees can be achieved during runtime by safety filtering for individual states \cite{hsunguyen2023isaacs, li2023learningpredictivesafetyfilter, wang2024magicsadversarialrlminimax, nguyen2024gameplay}. By designing a new discounted RA value function and a set-based certification method, \cite{li2024certifiable} addresses these issues.

On the other hand, there is no equivalent work for continuous-time systems that considers all of these settings. For example, using HJ-based methods, the authors in \cite{Zidani_2010} study the finite-horizon RA game (two-player). Although an extension to the infinite-horizon case is studied in \cite[Sec. 4]{Altarovici_Bokanowski_Zidani_2013}, it considers a running cost problem and imposes no constraints on the terminal state, so the resulting value function does not encode the information on whether the system can reach the target in some finite time, which is different from the objective in our work. We make a remark on the differences in more detail after we introduce our value function in Section III. The authors in \cite{Xue_ACC} compute the lower and upper robust controlled invariant sets for differential games, but only consider the avoidance problem. The verification of reach-avoid properties for autonomous systems is discussed in \cite{Xue_TAC}, but does not consider the synthesis of controllers for zero-sum games. The authors in \cite{Kene_2018} devise a contractive Bellman operator but only consider the reach or avoid case. Inspired by these prior works, we propose a new HJ-based method to solve the infinite-horizon two-player RA game that yields a contractive Bellman backup operator. By simultaneously accounting for all these problem components, which are often addressed separately in prior literature, we believe our work fills a large gap in this topic.

%have been proposed, where the RA set is characterized by the zero sublevel set of a value function which is the solution to a modified Hamilton-Jacobi variational inequality. Nonetheless, they focus on the cases ofMoreover, \cite{Altarovici_Bokanowski_Zidani_2013} relies on augmenting the state space by one more dimension, which introduces additional computational complexity that grows exponentially with problem dimensionality.

%For infinite-horizon games, the usage of a discounted formulation of the HJ-PDE is necessary for the existence of a unique solution \cite{Wabersich_Taylor_Choi_Sreenath_Tomlin_Ames_Zeilinger_2023}. 

%In the standard reach-avoid formulation, 
Although reach-avoid games provide a fundamental framework for ensuring a system reaches a target while avoiding unsafe states, they inherently focus on goal achievement rather than long-term safety; there are no guarantees on remaining in the goal (or maintaining safety) once the goal has been reached or after the prescribed time horizon. To ensure that the system remains in the goal set, the notion of reach-avoid-stay (RAS) has been extensively studied, and a popular family of methods tries to find a joint Control Lyapunov Barrier Function (CLBF), though each of them faces distinct challenges. For discrete-time systems, \cite{so2023solving, Chenevert_Li_Kannan_Bae_Lee_2024} resort to reinforcement learning for scalability, but neither provides verification methods for the learned RAS sets, i.e., safety and stay are not guaranteed in practice.

For continuous-time systems, \cite{ROMDLONY201639, MENG2022110478} lays the theoretical foundation for the existence of a CLBF if a given set satisfies the RAS condition. However, it assumes the system to be affine in both control and disturbance and does not provide a constructive way to find a CLBF, which is a highly difficult task for high-dimensional or complex systems due to the lack of a universal construction method. \cite{pmlr-v164-dawson22a} attempts to address this by leveraging the fitting power of neural networks to learn CLBFs, but the training process relies on the ability to sample from the control-invariant set. However, this becomes challenging when dealing with complex nonlinear dynamics, as the control-invariant set is not explicitly known. Note that for previous work on CLBFs, reach-avoid-stay is sometimes called interchangeably with stabilize-avoid (SA), as the Lyapunov property ensures stabilization of the system.
%but the validity of the Lyapunov property of the computed CLBF Lyapunov property ensures stabilization of the system.
%\cite{so2023solving} resorts to the epigraph form and reinforcement learning, but it does not consider disturbance. \cite{Chenevert_Li_Kannan_Bae_Lee_2024} takes disturbance into account
%Since the majority of the prior work relies on Control-Lyapunov-Barrier based methods, reach-avoid-stay is sometimes called interchangeably with stabilize-avoid (SA), as the Lyapunov property ensures stabilization of the system. For discrete-time systems, 

%\shnote{other work on stabilizing to goals safely (e.g. CLF/CBF methods, joint CLBF, anything else)}

In this paper, we propose a new HJ-based method to solve both the infinite-horizon RA and SA games for general nonlinear continuous-time systems. Specifically, our main contributions are:
%In this letter, we propose a new HJ-based method to solve both the reach-avoid and stabilize-avoid problem, under the settings of infinite-horizon and general nonlinear continuous-time systems subject to bounded disturbances. Specifically, our main contributions are:
\begin{enumerate}
    \item We define a new RA value function $V_\gamma$ and prove that the zero sublevel set of the RA value function is exactly the RA set.
    \item We establish the theoretical foundation to compute $V_\gamma$: the dynamic programming principle (DPP),  $V_\gamma$ being a unique viscosity solution to a Hamilton-Jacobi-Isaac's Variational Inequality (HJI-VI), and the contraction property of the Bellman backup operator associated with the DPP of $V_\gamma$.
    \item We present a two-step framework to construct the SA value function, and show that the zero sublevel set of the SA value function fully recovers the desired SA set.
    \item We show how to synthesize the controllers for the RA and SA tasks.
\end{enumerate}
To the best of our knowledge, both the proposed RA and SA frameworks are the first of their kind to solve the respective infinite-horizon games for general nonlinear continuous-time systems.
%problems under these settings.
%Our framework involves constructing a novel discounted reach-avoid method under these settings, which is then applied to a two-step method with the Robust Control Lyapunov Value Function (R-CLVF) \cite{gong2024robustcontrollyapunovvaluefunctions}. To the best of our knowledge, both the discounted reach-avoid and the stabilize-avoid frameworks we proposed are the first of their kind to solve the respective games under these settings.

The rest of this article is organized as follows. In Section \ref{sec: background}, we define the RA and SA problems and introduce the concept of Robust Control Lyapunov-Value Functions (R-CLVF) for stabilization. In Section \ref{sec: RA}, we introduce the value function $V_\gamma$ to solve the RA game, and discuss its HJ characterizations. In Section \ref{sec: SA}, we integrate our RA formulation with R-CLVF and propose a two-step framework to solve the SA game.
%\shnote{make contributions clear, consider bullet points}

\section{BACKGROUND}\label{sec: background}
In this section, we first define the system of interest alongside the RA and SA problems. 

\subsection{Dynamic System}
We consider a general nonlinear time-invariant system
\begin{align}\label{eq: dynamics}
    \dot x(t)  = f (x(t), \csig(t), \dsig(t)), \hspace{3mm} \text{for } t\geq0,\hspace{3mm} x(0)=x_0.
\end{align}
%where $t$ is the initial time and $x_0\in \R^n$ is the initial state,
The control signal \(\csig\colon[0, \infty) \rightarrow \cset\) and the disturbance signal \(\dsig\colon[0, \infty) \rightarrow \dset\) are drawn from the set of Lebesgue measurable control signals: \(
\cfset\coloneqq \{\csig\colon[0, \infty) \rightarrow \cset \mid \csig \text { is Lebesgue measurable.}\}\) and the set of Lebesgue measurable disturbance signals: \(
\dfset\coloneqq \{\dsig\colon[0, \infty) \rightarrow \dset \mid \dsig \text { is Lebesgue measurable.}\}\), where \(\cset \subset \mathbb{R}^{m_{u}}\) and \(\dset \subset \mathbb{R}^{m_{p}}\) are compact convex.

Assume the dynamics $f\colon\mathbb{R}^{n} \times \mathcal{\cset} \times \mathcal{D} \mapsto \R^{n}$ is uniformly continuous in \((x, u, d)\), Lipschitz continuous in $x, \forall u \in \cset, d \in \dset$, and bounded \(\forall x \in \mathbb{R}^{n}, u \in \cset, d \in \dset\). Under these assumptions, given the initial state \(x\), control and disturbance signal \(\csig, \dsig\), there exists a unique solution $\xi^{\csig, \dsig}_{x}(t)$ of the system \eqref{eq: dynamics}. When the initial condition, control, and disturbance signal used are not important, we use \(\xi(t)\) to denote the solution, which is also called the trajectory. Further assume the disturbance signal can be determined as a strategy with respect to the control signal: \(\lambda\colon \mathbb{U} \mapsto \mathbb{D}\), drawn from the set of non-anticipative maps \(\lambda \in \Lambda\).

\begin{definition} \(\Lambda\) The set of non-anaticipative maps $\Dmap$:
\begin{align*}
    \Dmap := \{ & \mathcal \lambda: \cfset\mapsto \dfset: \csig(\tvar) = \hat \csig(\tvar) \text{ a.e. } \forall \tvar \in [0,\infty) \\
    & \implies \mathcal \lambda[\csig](\tvar) = \mathcal \lambda[\hat \csig](\tvar) \text{ a.e. } \forall \tvar \in [0,\infty) \} .
\end{align*}
\end{definition}
With this non-anticipative assumption, the disturbance has an instantaneous advantage over the control signal.

\subsection{Reach-Avoid and Stabilize-Avoid Problems}
% In this section, we first introduce some necessary notations and then define the reach-avoid and stabilize-avoid problems we are trying to solve in this paper.

Let the open sets $\target \coloneqq \{x\colon \ell(x) < 0 \}\subseteq \R^n$ and $\constraint \coloneqq \{x\colon c(x) < 0 \}\subseteq \R^n$ denote the target set and constraint set, which are characterized by some Lipschitz continuous, bounded cost function $\ell(x)$ and constraint function $c(x)$, respectively. In this paper, we solve the following problems:

\textbf{Problem 1}: Specify a set of states that can be controlled to the target set safely in finite time (i.e., exists some finite $T\geq 0$ s.t. $\xi_{x}^{\csig, \lambda[\csig]}(T) \in \mathcal{T}$), under the worst-case disturbance. We call this set the \textit{\textbf{RA set}}:
%We are interested in finding two setsall the states that can be controlled to the target set safely in finite time, under \textit{the worst-case disturbance}. We call this set the \textit{\textbf{reach-avoid (RA) set}}:
\begin{align*}
    \mathcal{R} \mathcal{A}(\mathcal{T}, \mathcal{C})\coloneqq& \bigl\{x \in \mathbb{R}^{n}\colon \forall \lambda \in \Lambda, \exists \csig\in \cfset \text { and } T \geq 0, \\
    \text { s.t., } & \forall t\in [0, T], \xi_{x}^{\csig, \lambda[\csig]}(t) \in \constraint \text { and } \xi_{x}^{\csig, \lambda[\csig]}(T) \in \mathcal{T}\bigr\}.
\end{align*}
Note that though an individual state has to reach the target while satisfying constraints in a finite time, we do not require the existence of a prescribed finite $T$ that works for all states in the RA set. 

\textbf{Problem 2}: Specify a set of states that can be stabilized safely to a control invariant set in $\mathcal T$, under the worst-case disturbance. We call this set the \textit{\textbf{SA set}}:
%We are interested in finding two setsall the states that can be controlled to the target set safely in finite time, under \textit{the worst-case disturbance}. We call this set the \textit{\textbf{reach-avoid (RA) set}}:
\begin{align*}
    \mathcal{S} \mathcal{A}(\mathcal{T}, \mathcal{C})\coloneqq& \bigl\{x \in \mathbb{R}^{n}\colon \forall \lambda \in \Lambda, \exists \csig\in \cfset \text { and } T \geq 0,\\  \text { s.t., } & \forall t\geq 0,
    \xi_{x}^{\csig, \lambda[\csig]}(t) \in \constraint \text { and } \xi_{x}^{\csig, \lambda[\csig]}(T)\in \mathcal{T},\\ & \lim_{t\rightarrow \infty }\min_{y\in \partial \mathcal I_\text{m}} \| \xi_{x}^{\csig,  \lambda[\csig]}(t) -y \| =0 \bigr\},
\end{align*}
where $\mathcal I_\text{m}$ is the smallest robustly control invariant set defined below, after we introduce the R-CLVF. The finite-horizon RA problem can be solved by HJ reachability if a finite time horizon $T$ is specified, while the infinite-horizon RA and SA problems remain open.
\begin{remark} \label{remark: motivation}
    We highlight the fundamental difference between the RA (Problem 1) and SA (Problem 2) objectives. The RA problem requires the system to reach $\mathcal{T}$ at some time $T < \infty$ while remaining safe in $\mathcal{C}$ only on the interval $[0, T]$. No guarantees are provided after $T$; the trajectory may subsequently leave $\mathcal{T}$ or $\mathcal{C}$. On the other hand, the SA problem requires the system to remain safe in $\mathcal{C}$ for all time ($t \ge 0$) and after reaching $T$, asymptotically stabilize to a robust control invariant subset $\mathcal{I}_m \subseteq \mathcal{T}$.
\end{remark}
To solve the SA problem, we will apply the R-CLVF~\cite{gong2024robustcontrollyapunovvaluefunctions}:
\begin{definition}\label{def: R_CLVF}\textbf{R-CLVF} $\clvf: \mathcal D_\theta \mapsto \mathbb R $ of \eqref{eq: dynamics} is
\begin{align*} %\label{eq:R-CLVF}
    \clvf (x;p) = \lim _{\tinit \rightarrow \infty} \sup _{\lambda \in \Lambda} \inf _{\csig\in \cfset} \{ \sup_{s\in [0, t]} e^{\theta s}  r( \xi(s) ; p ) \}.
\end{align*}
Here, $ \mathcal D_\theta \subseteq \mathbb R^n$ is the domain, $\theta \geq 0$ is a user-specified parameter that represents the desired decay rate, $r (x;p) = \|x-p\| -a$, where $p$ is the desired point that we want to stabilize to, and $a$ is a parameter that depends on the system dynamics and will be explained later. 
\end{definition}
When $\theta = 0$ and $r (x;p) = \| x -p \|$, the value $\clvf(x)$ represents the largest deviation from $p$ along the trajectory starting from $x$. The level set corresponding to the smallest value $\minval = \min_x \clvf(x;p)$ is the \textbf{smallest} robustly control invariant set (SRCIS) $\mathcal{I_{\text{m}}}$ of $p$, defined in~\cite{gong2024robustcontrollyapunovvaluefunctions}. 
% Here, `smallest' means the smallest value.  

When $\theta>0$ and $r (x;p) = \| x -p \| - \minval$ (i.e., $ a = \minval$), the R-CLVF value $\clvf(x)$ captures the largest exponentially amplified deviation of a trajectory starting from $x$ to the $\mathcal I_\text{m}$, under worst-case disturbance. If this value is finite, it means $x$ can be exponentially stabilized to $\mathcal{I_{\text{m}}}$ (Lem.~7 of~\cite{gong2024robustcontrollyapunovvaluefunctions}). 

\begin{theorem} \label{thrm: clvf_exp_decay}
    The relative state can be exponentially stabilized to the $\mathcal I_{\text{m}}$ from $\mathcal D_\theta \setminus \mathcal I_{\text{m}}$, if the R-CLVF exists in $\mathcal D_\theta$, i.e., $\exists k > 0$, $\forall t \geq 0$ 
    \begin{align} \label{eqn:exp_decay}
        \min _ {a \in \partial\mathcal I_\text{m}}\|\xi(t)-a\|  \leq ke^{-\theta t} \min _ {a \in \partial \mathcal I_\text{m}}\|x-a\|.
    \end{align}
\end{theorem}
For conciseness, we simplify $r(x;p)$ and$\clvf(x;p)$ to $r(x)$ and $\clvf(x)$, as $p$ is a hyperparameter. The R-CLVF can be computed by solving the following R-CLVF-VI until convergence
\begin{align*}
     0 = &\max \{ r(x) - \clvf(x), \notag \\ &  \max_{\dstb \in \dset}\min_{u\in \cset}  D_x\clvf(x) \cdot f(x, u, d) + \theta \clvf(x) \}.
\end{align*}
The R-CLVF optimal controller is 
\begin{align} \label{eqn:opt_ctrl}
    \pi_{H} = \argmax_{\dstb \in \dset}\min_{\csig \in \cset} D_x \clvf(x) \cdot f(x, u, d).
\end{align}
The benefit of R-CLVF is two-fold: 1) for each positive level set, it is guaranteed that $$ \max_{\dstb \in \dset}\min_{u\in \cset}   D_x\clvf(x) \cdot f(x, u, d) \leq - \theta \clvf(x) <0,$$ and therefore each positive level set is attractive. 2) For any $\theta \geq 0$, the zero sublevel set of the R-CLVF is the \textbf{largest} robustly control invariant \textbf{subset} of the zero sub-level set of $r(x)$.

\section{A NEW DISCOUNTED RA VALUE FUNCTION}\label{sec: RA}
In this section, we propose a discounted RA value function to solve the RA and SA problems, i.e., characterizing $\mathcal{R} \mathcal{A}(\mathcal{T}, \mathcal{C})$ and $\mathcal{S} \mathcal{A}(\mathcal{T}, \mathcal{C})$. The introduction of the discount factor is necessary to ensure the value function is the unique viscosity solution of the corresponding HJI-VI, and many other desirable properties, such as the Lipschitz continuity and the contraction of the associated Bellman backup operator, also arise from the discount factor.
%its properties to characterize $\mathcal{R} \mathcal{A}(\mathcal{T}, \mathcal{C})$.
\subsection{Definition and Properties}
\begin{definition}
A time-discounted RA value function $V_{\gamma}(x): \R^n \mapsto \R$ is defined as
    \begin{align}\label{eq: def-inf}
        V_{\gamma}(x) &\coloneqq \sup _{\lambda \in \Lambda}\inf_{\csig\in \cfset}\inf_{t\in[0, \infty)} \notag\\
        &\max\bigl\{e^{-\gamma t} \ell(\xi(t)), \max_{s\in [0, t]}e^{-\gamma s} c(\xi(s))\bigr\},
    \end{align}    
\end{definition}
\noindent where $\xi$ solves \eqref{eq: dynamics} with initial state $x$ and $\gamma >0$. Intuitively, the control seeks to decrease the value along the trajectory, while the disturbance seeks to increase it. The infimum over time captures whether the trajectory ever reaches the target set $\target$ and the maximum checks if it ever violates the safety constraints. That is, $V_{\gamma}(x)< 0$ implies that there exists a control signal $\csig$ that can steer the system from initial state $x$ to the target without violating the constraints, even under worst-case disturbance. Before formally proving this claim, we make the following remarks to better contextualize our work in the literature.
%\zgnote{ What is the domain of this function? and maybe add a proposition discussing the continuity?}

\begin{remark}(Comparison with the formulation in \cite{Altarovici_Bokanowski_Zidani_2013}) \label{remark: literature}
    We want to point out that the infinite-horizon RA game and the value function we considered and constructed here in \textbf{Problem 1} and \eqref{eq: def-inf} are different from those in \cite{Altarovici_Bokanowski_Zidani_2013}. Specifically, although both of us and \cite[Sec. 4 \& 5]{Altarovici_Bokanowski_Zidani_2013} consider avoidance or state constraints, we focus on whether the system can reach the target in finite time (a constraint on the terminal state), so we design a terminal cost rather than a running cost. On the other hand, they emphasize running costs and do not impose requirements on the terminal state. As a result, their computed infinite-horizon value function cannot inform whether the system can be steered into the target in finite time. This difference becomes more obvious if one compares \eqref{eq: def-inf} with their value function (cf. Eq. (4.2) and (4.3) therein), which are as follows if written in our notation,
        \begin{align*}
        \tilde{V}_\gamma(x) &= \min_{\csig\in \cfset} \biggl ( \int_{0}^\infty e^{-\gamma s} \ell(\xi(s), \csig (s)) ds \\
        &\qquad \qquad \qquad \big| \hspace{1em} \xi(s) \in \mathcal C \quad \forall s\in (0,\infty) \biggr), \\
        w(x,z) &= \max\biggl\{\min_{\csig\in \cfset} \biggl ( \int_{0}^\infty e^{-\gamma s} \ell(\xi(s), \csig (s)) ds -z  \biggr ),\\
        &\qquad \qquad \quad \sup_{s \in (0,\infty)} \biggl ( e^{-\gamma s} c(\xi(s)) \biggr)\biggr\}.
    \end{align*}
\end{remark}

\begin{remark}(Exact Penalization) \label{remark: exact exact penalization}
    In optimal control, the idea of encoding state constraints through a \emph{trajectory-wise maximum} (``max-cost'', $L^\infty$-type) term as in \eqref{eq: def-inf} is a classical approach \cite{BARRON19891067, Barron1999, Barron_Jensen_1996}. We adopt this method in our work as other Hamilton-Jacobi frameworks for state-constrained control developments \cite{Zidani_2010, Margellos_2011, Zidani_2011, Altarovici_Bokanowski_Zidani_2013}.   
\end{remark}

Now, we prove that the zero sublevel set of the value function \eqref{eq: def-inf} gives the $\mathcal{R} \mathcal{A}(\mathcal{T}, \mathcal{C})$ set we are looking for in \textbf{Problem 1}.
\begin{proposition} (Exact Recovery of RA set) \label{prop: def-RA}
    \begin{equation}\label{eq: def-RA}
        \mathcal{R} \mathcal{A}(\mathcal{T}, \mathcal{C}) = \{x\colon V_{\gamma}(x) < 0 \}    
    \end{equation}
    % \bynote{$\leq$ or $<$? For the current contrapositive statement, we need $\leq$.}
    % \zgnote{if you define $\mathcal C = x: c(x)<0$, then it should be $<$. }
    
    \begin{proof}
        We first prove the sufficiency. Suppose $x\in \mathcal{R} \mathcal{A}(\mathcal{T}, \mathcal{C})$. That is, given arbitrary $\lambda \in \Lambda, \exists \hat{\csig}\in \cfset \text { and } T \geq 0$ such that $\max\limits_{t\in [0, T]} c(\xi^{\hat{\csig}, \lambda[\hat{\csig}]}_x(t)) < 0$ and $\ell(\xi^{\hat{\csig}, \lambda[\hat{\csig}]}_x(T))<0$. Since $e^{-\gamma t}>0, \forall t\in \R$, we have $\max\limits_{t\in [0, T]}e^{-\gamma t} c(\xi^{\hat{\csig}, \lambda[\hat{\csig}]}_x(t)) < 0$ and $e^{-\gamma T}\ell(\xi^{\hat{\csig}, \lambda[\hat{\csig}]}_x(T))<0$. Therefore, 
        \begin{align*}
            \inf\limits_{\csig\in \cfset}&\inf\limits_{T\in[0, \infty)} \notag\\
            &\max\{e^{-\gamma T}\ell(\xi^{\csig, \lambda[\csig]}_x(T)), \max\limits_{t\in [0, T]}e^{-\gamma t} c(\xi^{\csig, \lambda[\csig]}_x(t))\} < 0.
        \end{align*}
 Because $\lambda$ is arbitrary, $ V_{\gamma}(x) < 0$.
%        \begin{align*}
 %           V_{\gamma}(x) &\coloneqq \sup _{\lambda \in \Lambda}\inf_{\csig\in \cfset}\inf_{t\in[0, \infty]} \notag\\
 %           &\max\{e^{-\gamma t} \ell(\xi(t)), \sup_{s\in [0, t]}e^{-\gamma s} c(\xi(s))\}<0.
 %       \end{align*}

    %For the necessity, we prove the contrapositive statement. Let $x$ be such that $V_{\gamma}(x) \geq 0$, i.e., $\exists\lambda \in \Lambda$, such that $\forall \csig\in \cfset \text { and } T \geq 0$, 
    % \zgnote{this statement is not rigorous: $\sup_ \lambda (something) \geq 0$ does not mean $\exists \lambda $ s.t. something $\geq 0$. eg. $\sup _{x\geq 0} -\frac{1}{x} = 0$, but exists no $x\geq 0$ s.t. $-\frac{1}{x} = 0$. This is why we need $\leq$ in (8) right?}
    For the necessity, let $x$ be such that $V_{\gamma}(x) < 0$, that is, for any $\lambda \in \Lambda, \exists \hat{\csig}\in \cfset$ such that, 
    \begin{align}\label{eq: prop 1 necessity}
        \inf_{t\in[0, \infty)}\max\{e^{-\gamma t}\ell(\xi^{\hat{\csig}, \lambda[\hat{\csig}]}_x(t)),
        \max\limits_{s\in [0, t]}e^{-\gamma t} c(\xi^{\hat{\csig},\lambda[\hat{\csig}]}_x(s))\} < 0.
    \end{align}
    Thus, there also exists $T\geq 0$ such that
    \begin{equation*}
    \max\{e^{-\gamma T}\ell(\xi^{\hat{\csig}, \lambda[\hat{\csig}]}_x(T)), 
        \max\limits_{s\in [0, T]}e^{-\gamma t} c(\xi^{\hat{\csig},\lambda[\hat{\csig}]}_x(s))\} < 0.        
    \end{equation*}
     Since $e^{-\gamma t}>0$, both $\ell(\xi^{\hat{\csig}, \lambda[\hat{\csig}]}_x(T)) $ and $\max\limits_{t\in [0, T]} c(\xi^{\hat{\csig}, \lambda[\hat{\csig}]}_x(t)) $ are negative, which means $x\in \mathcal{R} \mathcal{A}(\mathcal{T}, \mathcal{C})$.
    \end{proof}
\end{proposition}
%In fact, many desirable properties of the value function, such as its Lipschitz continuity and the contraction of the associated Bellman backup operator, arise from the discount factor.
We want to point out that for this proof, the \textit{openness} in the definition of $\target$ and $C$ is the key to proving the necessity part. Specifically, from \eqref{eq: prop 1 necessity}
%\[
%\inf_{t\in[0, \infty)}\max \bigl\{ e^{-\gamma T} \ell(\xi(T)),\; \max_{s \in [0,T]} e^{-\gamma s} c(\xi(s)) \bigr\} < 0
%\]
we concluded that $\ell(\xi(T)) < 0$ and $c(\xi(s)) < 0$ for all $s \in [0,T]$, hence
$\xi(T) \in \target$ and $\xi(s) \in \constraint$ for all $s \in [0,T]$. If, instead, $\target$ and $\constraint$ were taken to be \textit{closed} sets, the same argument would make the strict inequality in \eqref{eq: prop 1 necessity} be non-strict, which does not exclude the case $\ell(\xi(T)) > 0$ or $c(\xi(s)) > 0$ and therefore does not allow us to assert $x \in  \mathcal{R} \mathcal{A}(\mathcal{T}, \mathcal{C})$. This is why we explicitly assume $\target$ and $\constraint$ to be open in their definitions.
\begin{remark} [Identical zero sublevel set.]\label{remark:same_zero_set}
    Note that in this proof, we do not make any assumption on $\gamma$ other than being positive. In other words, Proposition~\ref{prop: def-RA} holds for any $\gamma > 0$, i.e., the RA value functions~\eqref{eq: def-inf} with different $\gamma$'s share the same zero sublevel set, which is exactly the RA set.
\end{remark}

As mentioned before, we now show how the discount factor guarantees the boundedness and Lipschitz continuity of~\eqref{eq: def-inf}.
\begin{proposition}(Boundedness and Lipschitz Continuity)\label{prop: V_RA properties}
    $V_{\gamma}$ is bounded and Lipschitz continuous in $\R^n$ if $L_f < \gamma$, where $L_f$ is the Lipschitz constant of $f$.
    
    \begin{proof}
        Define
        \begin{equation*}
            p(\lambda, \csig, x, t) = \max\bigl\{e^{-\gamma t} \ell(\xi(t)), \max_{s\in [0, t]}e^{-\gamma s} c(\xi(s))\bigr\}.
        \end{equation*}
        Take $x_1, x_2\in \R^n$ and any $ \varepsilon > 0$. Then, there exists $\hat{\dmap}\in \Dmap$ such that for any $\csig\in \cfset$,
        \begin{equation*}
            V_{\gamma}(x_1) \leq \inf_{t\in[0, \infty)} p(\hat{\dmap}, \csig, x_1, t) + \varepsilon.
        \end{equation*}
        Similarly, for any $\dmap \in \Dmap$, there exists $\hat{\csig}\in \cfset$ such that
        \begin{equation*}
            V_{\gamma}(x_2) \geq \inf_{t\in[0, \infty)} p(\dmap, \hat{\csig}, x_2, t) - \varepsilon.
        \end{equation*}
        Combining the above two inequalities, we have
        \begin{align}\label{eq: Lip_1}
            V_{\gamma}(x_1) - V_{\gamma}(x_2) \leq &\inf_{t\in[0, \infty)} p(\hat{\dmap}, \hat{\csig}, x_1, t) \notag\\
            & \quad - \inf_{t\in[0, \infty)} p(\hat{\dmap}, \hat{\csig}, x_2, t) + 2\varepsilon.    
        \end{align}
        Moreover, there exists $\hat{t}\in [0, \infty)$ such that
        \begin{align*}
            &\inf_{t\in[0, \infty)} p(\hat{\dmap}, \hat{\csig}, x_2, t) \geq p(\hat{\dmap}, \hat{\csig}, x_2, \hat{t}) - \varepsilon,\\
            &\inf_{t\in[0, \infty)} p(\hat{\dmap}, \hat{\csig}, x_1, t) \leq p(\hat{\dmap}, \hat{\csig}, x_1, \hat{t}) .
        \end{align*}
        Plugging this back to \eqref{eq: Lip_1}, we have
        \begin{align}\label{eq: Lip_2}
            V_{\gamma}(x_1) - V_{\gamma}(x_2) \leq &\inf_{t\in[0, \infty)} p(\hat{\dmap}, \hat{\csig}, x_1, t) \notag\\
            & \quad - p(\hat{\dmap}, \hat{\csig}, x_2, \hat{t}) + 3\varepsilon, \notag\\
            \leq &p(\hat{\dmap}, \hat{\csig}, x_1, \hat{t}) - p(\hat{\dmap}, \hat{\csig}, x_2, \hat{t}) + 3\varepsilon.
        \end{align}
        By definition of $p$, both of the following hold:
        \begin{align}
            p(\hat{\dmap}, \hat{\csig}, x_2, \hat{t}) &\geq e^{-\gamma \hat{t}} \ell(\xi^{\hat{\csig}, \hat{\dmap}[\hat{\csig}]}_{x_2}(\hat{t})) ,\label{eq: Lip p2_l}\\
            p(\hat{\dmap}, \hat{\csig}, x_2, \hat{t}) &\geq \max_{s\in [0, \hat{t}]}e^{-\gamma s} c(\xi^{\hat{\csig}, \hat{\dmap}[\hat{\csig}]}_{x_2}(s)), \notag\\
            & \geq e^{-\gamma s} c(\xi^{\hat{\csig}, \hat{\dmap}[\hat{\csig}]}_{x_2}(s)), \label{eq: Lip p2_c}
        \end{align}
        for any $s\in [0, \hat{t}]$. Similarly, either of the following holds:
        \begin{align}
            p(\hat{\dmap}, \hat{\csig}, x_1, \hat{t}) &\leq e^{-\gamma \hat{t}} \ell(\xi^{\hat{\csig}, \hat{\dmap}[\hat{\csig}]}_{x_1}(\hat{t})) ,\label{eq: Lip p1_l}\\
            p(\hat{\dmap}, \hat{\csig}, x_1, \hat{t}) &\leq \max_{s\in [0, \hat{t}]}e^{-\gamma s} c(\xi^{\hat{\csig}, \hat{\dmap}[\hat{\csig}]}_{x_1}(s)), \notag\\
            & \leq e^{-\gamma \hat{s}} c(\xi^{\hat{\csig}, \hat{\dmap}[\hat{\csig}]}_{x_1}(\hat{s})), \label{eq: Lip p1_c}
        \end{align}
        for some $\hat{s}\in [0, \hat{t}]$. Finally, plugging \eqref{eq: Lip p2_l} and \eqref{eq: Lip p1_l} back to \eqref{eq: Lip_2}, we have
        \begin{align}\label{eq: Lip_3}
            V_{\gamma}(x_1) - V_{\gamma}(x_2) &\leq  e^{-\gamma \hat{t}} \ell(\xi^{\hat{\csig}, \hat{\dmap}[\hat{\csig}]}_{x_1}(\hat{t})) \notag\\
            & \qquad - e^{-\gamma \hat{t}} \ell(\xi^{\hat{\csig}, \hat{\dmap}[\hat{\csig}]}_{x_2}(\hat{t})) + 3\varepsilon, \notag\\
            &\leq L_{\ell}e^{-\gamma \hat{t}}e^{L_f\hat{t}}||x_1 - x_2|| + 3\varepsilon, \notag\\
            &\leq L_{\ell}||x_1 - x_2|| + 3\varepsilon,             
        \end{align}
        where $L_{\ell}, L_f$ are the Lipschitz constants of $\ell$ and $f$, respectively. The second inequality is a result of Gronwall's inequality, and the third is due to the assumption that $L_f < \gamma$. Similarly, we can show that $ V_{\gamma}(x_2) - V_{\gamma}(x_1) \leq L_{\ell}||x_1 - x_2|| + 3\varepsilon$. Thus, given any two states $x_1$ and $x_2$, since $\epsilon$ can be choosen arbitrarily small, we have $|V_{\gamma}(x_1) - V_{\gamma}(x_2)| \leq L_{\ell}||x_1 - x_2||$. The proof for plugging \eqref{eq: Lip p2_c} and \eqref{eq: Lip p1_c} to \eqref{eq: Lip_2} is similar.
        
        The boundedness of $V_\gamma$ follows from that of $\ell(x)$ and $c(x)$, and the fact that $\lim_{t\to \infty} e^{-\gamma t} = 0$.
    \end{proof}
\end{proposition}
Note that since the value function is Lipschitz continuous, it is differentiable almost everywhere by Rademacher’s Theorem \cite[Ch.5.8.3]{evans2010chapter9.2}.
\vspace{-1em}
\subsection[Hamilton-Jacobi Characterization of V gamma]{Hamilton-Jacobi Characertization of $V_\gamma$}
We now establish the theoretical foundations to compute $V_{\gamma}$: the DPP as a result of Bellman's optimality principle, $V_{\gamma}$ being a unique viscosity solution to a HJI-VI, and the contraction property of the Bellman backup operator associated with the DPP of $V_{\gamma}$. These provide two approaches for the numerical computation. The first is to follow the procedure in~\cite{Fisac_Chen_Tomlin_Sastry_2015}, i.e., solving the HJI-VI and combining the DPP until convergence. The other is to do a value iteration based on the Bellman backup operator. In this paper, all numerical solutions are obtained following the first approach. 

We first show the DPP, as it is the basis for proving the viscosity solution and deriving the Bellman backup operator. 
\begin{theorem}\label{th: dpp}
    (Dynamic Programming Principle). Suppose \(\gamma>0\). For \(x \in \mathbb{R}^{n} \text{ and } T>0\),
    \begin{align}\label{eq: DPP}
        V_{\gamma}(x)&=\sup _{\lambda \in \Lambda}\inf_{\csig\in \cfset} \min \bigl\{ \notag\\
        &\min _{t\in [0, T]} \max [e^{-\gamma t}\ell(\xi(t)), \max _{s \in[0, t]} e^{-\gamma s} c(\xi(s))], \notag\\
        &\max[e^{-\gamma T} V_{\gamma}(\xi(T)), \max _{t\in [0, T]} e^{-\gamma t} c(\xi(t))]\bigr\}.
    \end{align}
    \begin{proof}
    %\vspace{-1em}
        Separating the infimum from $[0,\infty)$ to the minimum between $\min_{t\in[0,T]}$ and  $\inf_{t\in[T,\infty)}$, we rewrite $V_{\gamma}(x)$ as,
        \begin{align}\label{eq: DPP_proof}
            V_{\gamma}(x)&=\sup _{\lambda \in \Lambda}\inf_{\csig\in \cfset} \min \bigl\{ \notag\\
            &\min_{t\in [0, T]} \max [e^{-\gamma t}\ell(\xi(t)), \max_{s \in[0, t]} e^{-\gamma s} c(\xi(s))], \notag\\
            &\inf_{t \in[T, \infty)} \max [e^{-\gamma t}\ell(\xi(t)),
            %\textcolor{red}{\sup_{s \in[T, t] ?}}
            \max_{s \in[0, t]}
            e^{-\gamma s} c(\xi(s))]\bigr\}.
        \end{align}
    Thus, it suffices to prove that the optimized last term (for both control and disturbance) in the minimum of \eqref{eq: DPP_proof} is identical to that of \eqref{eq: DPP}. Due to the time-invariance of the system \eqref{eq: dynamics}, we can transform
    \begin{align*}
        &\inf_{t \in[T, \infty)} \max [e^{-\gamma t}\ell(\xi(t)), \max_{s \in[0, t]} e^{-\gamma s} c(\xi(s))]\notag \\
        =& \inf_{t \in[T, \infty)} \max [e^{-\gamma t}\ell(\xi(t)), \max_{s \in[T, t]} e^{-\gamma s} c(\xi(s)),\notag \\
        & \max _{s \in[0, T]} e^{-\gamma s} c(\xi(s)) ],\notag \\
        =& \inf_{t \in[T, \infty)} \max \bigl \{ \max [ e^{-\gamma t}\ell(\xi(t)), \max _{s \in[T, t]} e^{-\gamma s} c(\xi(s)) ],  \notag \\
        & \max_{s \in[0, T]} e^{-\gamma s} c(\xi(s))  \bigr \}. 
\end{align*}
For the first equation, we use the same trick on separating the interval $[0,t]$ into $[0,T]$ and $[T,t]$. The second equation follows from extracting and adding a maximum for the comparison between the three terms.

Taking $k\coloneqq t-T$ and $y\coloneqq \xi^{\csig, \dsig}_{x}(T)$, we get 
\begin{align*}
         & \inf_{k \in[0, \infty)} \max \bigl \{ \max [ e^{-\gamma (k+T)}\ell(\xi^{\csig, \dsig}_{y}(k)),   \notag \\
         & \max _{s \in[0, k]} e^{-\gamma (s+T)} c(\xi^{\csig, \dsig}_{y}(k)) ], \max _{s \in[0, T]} e^{-\gamma s} c(\xi(s))  \bigr \}.   \\
        %& \inf_{t \in[T, \infty)} \max \bigl \{ \max [ e^{-\gamma t}\ell(\xi(t)), \max _{s \in[T, t]} e^{-\gamma s} c(\xi(s)) ],  \notag \\
        %& \max_{s \in[0, T]} e^{-\gamma s} c(\xi(s))  \bigr \}\\
        =&   \max \bigl \{ \inf_{k \in[0, \infty)} \max [ e^{-\gamma (k+T)}\ell(\xi^{\csig, \dsig}_{y}(k)), \notag \\
        &\max _{s \in[0, k]} e^{-\gamma (s+T)} c(\xi^{\csig, \dsig}_{y}(k)) ], \max _{s \in[0, T]} e^{-\gamma s} c(\xi(s))  \bigr \}, \notag 
        % &=\max\{\inf_{t \in[T, \infty]} [e^{-\gamma t}\ell(\xi(t)),           \notag \\ 
        % &\qquad \sup_{s \in[T, t]} e^{-\gamma s} c(\xi(s))], \sup_{s \in[0, T]} e^{-\gamma s} c(\xi(s))\} \notag \\ 
        % &= \max \{\inf_{k \in[0, \infty]} e^{-\gamma T}[e^{-\gamma k}\ell(\xi^{\csig, \dsig}_{y}(k)),\notag \\ 
        % &\qquad\sup_{s \in[0, k]} e^{-\gamma s} c(\xi^{\csig, \dsig}_{y}(s))], \
        % \sup_{s \in[0, T]} e^{-\gamma s} c(\xi(s))\}. \notag          
     %   &= \max\bigl\{\sup _{\lambda \in \Lambda}\inf_{\csig\in \cfset}\inf_{k \in[0, \infty]} \max e^{-\gamma T}[e^{-\gamma k}\ell(\xi^{\csig, \dsig}_{y}(k)),\notag \\ 
     %   &\qquad\sup_{s \in[0, k]} e^{-\gamma s} c(\xi^{\csig, \dsig}_{y}(s))], \
      %   \sup _{\lambda \in \Lambda}\inf_{\csig\in \cfset}\sup_{s \in[0, T]} e^{-\gamma s} c(\xi(s))\} \notag \\
     %   &= \max\bigl\{e^{-\gamma T}V_{\gamma}(\xi^{\csig, \dsig}_{x}(T)), \ \sup _{\lambda \in \Lambda}\inf_{\csig\in \cfset}\sup_{s \in[0, T]} e^{-\gamma s} c(\xi(s))\} 
    \end{align*}
    where the infimum over $k$ for $ \max _{s \in[0, T]} e^{-\gamma s} c(\xi(s))$ is dropped because it is independent of $k$. Plugging the above into \eqref{eq: DPP_proof} and noticing that the optimization over the two intervals $t \in [0,T]$ and $t \in [T,\infty) $ is independent, we have
    \allowdisplaybreaks
    \begin{align}
        V_{\gamma}(x)=&\sup _{\lambda \in \Lambda}\inf_{\csig\in \cfset} \min \biggl\{ \notag\\
        &\min_{t\in [0, T]} \max [e^{-\gamma t}\ell(\xi(t)), \max _{s \in[0, t]} e^{-\gamma s} c(\xi(s))], \notag\\
        &  \max \bigl \{ \inf_{k \in[0, \infty)} \max [ e^{-\gamma (k+T)}\ell(\xi^{\csig, \dsig}_{y}(k)), \notag \\
        & \max _{s \in[0, k]} e^{-\gamma (s+T)} c(\xi^{\csig, \dsig}_{y}(k)) ], \max _{s \in[0, T]} e^{-\gamma s} c(\xi(s))  \bigr \}  \biggr\} ,\notag \\
        % &\max \{\min _{k \in[0, \infty]} e^{-\gamma T}[e^{-\gamma k}\ell(\xi^{\csig, \dsig}_{y}(k)),\notag \\ 
        % &\qquad\sup_{s \in[0, k]} e^{-\gamma s} c(\xi^{\csig, \dsig}_{y}(s))], \
        % \sup_{s \in[0, T]} e^{-\gamma s} c(\xi(s))\}\bigr\} \notag \\
        =& \sup _{\lambda \in \Lambda}\inf_{\csig\in \cfset} \min \biggl\{  \notag\\
        & \min _{t\in [0, T]} \max [e^{-\gamma t}\ell(\xi(t)), \max _{s \in[0, t]} e^{-\gamma s} c(\xi(s))],  \notag\\
        &   \max \bigl \{ \sup _{\lambda \in \Lambda}\inf_{\csig\in \cfset} \inf_{k \in[0, \infty)} \max [ e^{-\gamma (k+T)}\ell(\xi^{\csig, \dsig}_{y}(k)), \notag \\
        & \max _{s \in[0, k]} e^{-\gamma (s+T)} c(\xi^{\csig, \dsig}_{y}(k)) ],  \max _{s \in[0, T]} e^{-\gamma s} c(\xi(s))  \bigr \}  \biggr\} ,\notag \\
        =&\sup _{\lambda \in \Lambda}\inf_{\csig\in \cfset} \min \biggl\{ \notag\\
        &\min _{t\in [0, T]} \max [e^{-\gamma t}\ell(\xi(t)), \max _{s \in[0, t]} e^{-\gamma s} c(\xi(s))], \notag\\
        &\max\{e^{-\gamma T} V_{\gamma}(\xi(T)), \max _{t\in [0, T]} e^{-\gamma t} c(\xi(t))\}\biggr\},        
    \end{align}
    where for the second equation, we plug in the definition of $V_\gamma(\xi(T))$. Notice the last equation is identical to the expression in \eqref{eq: DPP}.
   % since the first term in the maximum can be optimized independently of the rest, once $y\coloneqq \xi^{\csig, \dsig}_{x}(T)$ is specified.
    \end{proof}
\end{theorem} 
Building on Theorem \ref{th: dpp}, Theorem \ref{th: HJI} shows that~\eqref{eq: def-inf} is the unique viscosity solution to the HJI-VI defined below.  We follow the definition of the viscosity solution provided in~\cite{BARRON19891067}.
\begin{theorem}\label{th: HJI}(Unique Viscosity Solution to HJI-VI). Assume \(f\) satisfies \eqref{eq: dynamics}, and that \(\ell(x), c(x)\) are bounded and Lipschitz continuous. Then the value function \(V_{\gamma}(x)\) defined in \eqref{eq: def-inf} is the unique viscosity solution of the following HJI-VI
    \begin{align}\label{eq: HJI-VI}
        0=\max\bigl\{&\min \{\max _{d \in \dset}\min_{u\in \cset} D_xV_{\gamma}(x) \cdot f(x, u, d) - \gamma V_{\gamma}(x), \notag\\
        & \ell(x)-V_{\gamma}(x)\}, c(x)-V_{\gamma}(x)\bigr\}. 
    \end{align}
    
    \begin{proof}
        The structure of the proof follows the classical
        approach in \cite{Bardi_Capuzzo-Dolcetta_2009}, analogously to \cite{Fisac_Chen_Tomlin_Sastry_2015}.
        % , and draws from viscosity solution theory.
        
        A continuous function is a viscosity solution of a partial differential equation if it is both a subsolution and a supersolution (defined below). We will first prove that \(V_{ \gamma}\) is a viscosity subsolution of \eqref{eq: HJI-VI}. Let \(\psi \in C^{1}(\R^n)\) be such that \(V_{ \gamma}-\psi\) attains a local maximum at \(x_{0}\); without loss of generality, assume that this maximum is 0 . We say that \(V_{ \gamma}\) is a subsolution of \eqref{eq: HJI-VI} if, for any such \(\psi\),
        \begin{align}\label{eq: sub sol}
            \max\bigl\{&\min \{ H\left( x_{0}, D_{x} \psi(x_{0})\right) - \gamma \psi(x_0), \notag\\
            & \ell(x_0)-\psi(x_0)\}, c(x_0)-\psi(x_0)\bigr\}\geq 0,
        \end{align}
        where the Hamiltonian is defined as
        \begin{equation*}
            H\left( x, D_{x}\psi(x) \right)\coloneqq \max _{d \in \dset}\min_{u\in \cset} D_x\psi(x) \cdot f(x, u, d).
        \end{equation*}
%    Note that $H$ is continuous in $(x, \psi)$ as ensured by the compactness of $\cset$ and $\dset$ and the continuity of $f$.
    
    Assuming \eqref{eq: sub sol} is false, then the following must hold: 
    \begin{equation}\label{eq: c-psi}
         c(x_0)\leq \psi(x_0) - \varepsilon_1.
    \end{equation}
    Moreover, at least one of the following must hold:
   \begin{subequations}
        \begin{align}
            &\ell(x_0) \leq \psi(x_0) - \varepsilon_2 \label{eq: r-psi},\\
            &H\left( x_{0}, D_{x} \psi(x_{0})\right) - \gamma \psi(x_0)\leq - \varepsilon_3, \label{eq: H-psi}
        \end{align}
   \end{subequations}
   for some $\varepsilon_1, \varepsilon_2, \varepsilon_3>0$. We will show that these cannot be true, i.e., \eqref{eq: c-psi} and one of \eqref{eq: r-psi} and \eqref{eq: H-psi} cannot hold simultaneously. Suppose \eqref{eq: c-psi} and \eqref{eq: r-psi} are true. We abbreviate $\xi^{\csig, \dsig}_{x_0}(\cdot)$ to $\xi(\cdot)$ whenever the statement holds for any $\csig, \dsig$. By continuity of $\ell, c$ and system trajectories, there exists small enough $\delta > 0$, such that for all $\csig(\cdot), \dsig(\cdot), \tau\in [0, \delta]$,
   \begin{subequations}\label{eq: sub_c+l<=V}
    \begin{align}
       &e^{-\gamma \tau}c(\xi(\tau))\leq \psi(x_0) - \frac{\varepsilon_1}{2} = V_{ \gamma}(x_0) - \frac{\varepsilon_1}{2}, \label{eq: sub_c<=V}\\
       &e^{-\gamma \tau}\ell(\xi(\tau))\leq \ell(x_0) - \frac{\varepsilon_2}{2} = V_{ \gamma}(x_0) - \frac{\varepsilon_2}{2}.\label{eq: sub_l<=V}
   \end{align}
   \end{subequations}
   Incorporating \eqref{eq: sub_c+l<=V} into the dynamic programming principle \eqref{eq: DPP}, we have 
   \begin{align*}
       V_{ \gamma}(x_0)&\leq \sup _{\lambda \in \Lambda}\inf_{\csig\in \cfset} \min \bigl\{ \\
    &\min_{t\in [0, \delta]} \max [e^{-\gamma t}\ell(\xi(t)), \max_{s \in[0, t]} e^{-\gamma s} c(\xi(s))]\bigr\} \\
    &\leq V_{\gamma}(x_0) - \min \left\{\frac{\varepsilon_1}{2}, \frac{\varepsilon_2}{2}\right\},
   \end{align*}
   which is a contradiction, since $\varepsilon_1, \varepsilon_2>0$.
   
   Now suppose \eqref{eq: c-psi} and \eqref{eq: H-psi} are true. By definition of $H$, for any $\lambda\in \Lambda$
   \begin{align*}
       \inf_{\csig\in \cfset} D_x\psi(x_0) &\cdot f(x_0, \csig(0), \lambda[\csig](0)) - \gamma \psi(x_0),\\
       &\leq  H\left( x_{0}, D_{x} \psi(x_{0})\right) - \gamma \psi(x_0),\\
       &\leq -\varepsilon_3.
   \end{align*}
   Then, for small enough $\delta > 0$, \eqref{eq: sub_c<=V} holds and there exists $\bar{\csig}\in \cfset$ such that
   \begin{align*}
       D_x\psi(\xi_{x_0}^{\bar{\csig}, \lambda}(\tau)) &\cdot f(\xi_{x_0}^{\bar{\csig}, \lambda}(\tau), \bar{\csig}(\tau), \lambda[\bar{\csig}](\tau)) \\
       &- \gamma \psi(\xi_{x_0}^{\bar{\csig}, \lambda}(\tau))\leq  -\frac{\varepsilon_3}{2}, %\forall \tau\in [0, \delta].
   \end{align*}
   for all $\tau\in [0, \delta]$.
   Following the technique used in Theorem 1.10 in Chapter VIII in \cite{Bardi_Capuzzo-Dolcetta_2009} with some modifications, we multiply both sides of the above inequality by $e^{-\gamma \tau}$ and integrate from $0$ to $\delta$ to get
   \begin{align*}
       e^{-\gamma \delta}\psi(\xi_{x_0}^{\bar{\csig}, \lambda}(\delta)) - \psi(x_0) \leq -\frac{\varepsilon_3}{2}\delta.
   \end{align*}
   Recalling that $V_{\gamma}-\psi$ attains a local maximum of $0$ at $x_0$, we have
    \begin{align}\label{eq: sub_V<=eps}
       e^{-\gamma \delta}V_{\gamma}(\xi_{x_0}^{\bar{\csig}, \lambda}(\delta)) \leq  V_{\gamma}(x_0) - \frac{\varepsilon_3}{2}\delta.
   \end{align}
   Incorporating \eqref{eq: sub_c<=V} and \eqref{eq: sub_V<=eps} into the dynamic programming principle \eqref{eq: DPP}, we have 
   \begin{align*}
        V_{ \gamma}(x_0)&\leq \sup _{\lambda \in \Lambda}\inf_{\csig\in \cfset}
        \bigl\{\max[e^{-\gamma \delta} V_{\gamma}(\xi(\delta)), \max _{t\in [0, \delta]} e^{-\gamma t} c(\xi(t))]\bigr\},\\
        & \leq V_{\gamma}(x_0) - \min \left\{\frac{\varepsilon_1}{2}, \frac{\varepsilon_3}{2}\delta\right\}.
   \end{align*}
    which is a contradiction, since $\varepsilon_1, \varepsilon_3, \delta > 0$. Therefore, we conclude that \eqref{eq: sub sol} must be true and hence $V_{ \gamma}$ is indeed a subsolution of \eqref{eq: HJI-VI}.
    
    We now proceed to prove that \(V_{ \gamma}\) is a viscosity supersolution of \eqref{eq: HJI-VI}. Let \(\psi \in C^{1}(\R^n)\) be such that \(V_{ \gamma}-\psi\) attains a local minimum at \(x_{0}\); without loss of generality, assume that this minimum is 0. We say that \(V_{ \gamma}\) is a supersolution of \eqref{eq: HJI-VI} if, for any such \(\psi\),
        \begin{align}\label{eq: sup sol}
            \max\bigl\{&\min \{ H\left( x_{0}, D_{x} \psi(x_{0})\right) - \gamma \psi(x_0), \notag\\
            & \ell(x_0)-\psi(x_0)\}, c(x_0)-\psi(x_0)\bigr\}\leq 0,
        \end{align}
    Assume \eqref{eq: sup sol} is false, then either it must hold that
    \begin{equation}\label{eq: c+psi}
         c(x_0)\geq \psi(x_0) + \varepsilon_1,
    \end{equation}
    or both of the following are true:
   \begin{subequations}\label{eq:r&H+psi}
        \begin{align}
            &\ell(x_0) \geq \psi(x_0) + \varepsilon_2 \label{eq: r+psi},\\
            &H\left( x_{0}, D_{x} \psi(x_{0})\right) - \gamma \psi(x_0)\geq 2\varepsilon_3, \label{eq: H+psi}
        \end{align}
   \end{subequations}
   for some $\varepsilon_1, \varepsilon_2, \varepsilon_3>0$. Suppose \eqref{eq: c+psi} is true. Then, there exists small enough $\delta > 0$, such that for all $\csig(\cdot), \dsig(\cdot), \tau\in [0, \delta]$,
   \begin{equation}\label{eq: sup_c>=psi}
       e^{-\gamma \tau}c(\xi(\tau))\geq \psi(x_0) + \frac{\varepsilon_1}{2} = V_{ \gamma}(x_0) + \frac{\varepsilon_1}{2}.
   \end{equation}
   Incorporating \eqref{eq: sup_c>=psi} into the dynamic programming principle \eqref{eq: DPP}, we get 
    \begin{align*}
        V_{\gamma}(x) \geq \sup _{\lambda \in \Lambda}\inf_{\csig\in \cfset} \min \bigl\{&\min_{t\in [0, \delta]}\max_{s \in[0, t]} e^{-\gamma s} c(\xi(s)), \\
        &\max _{t\in [0, \delta]} e^{-\gamma t} c(\xi(t))\bigr\}, \\
        \geq V_{ \gamma}(x_0) + \frac{\varepsilon_1}{2}&,
    \end{align*}
    which is a contradiction, as $\varepsilon_1 > 0$. Now, suppose \eqref{eq:r&H+psi} holds. Then, there exists small enough $\delta_1 > 0$, such that for all $\csig(\cdot), \dsig(\cdot), \tau\in [0, \delta_1]$,
   \begin{equation}\label{eq:sup_l>=V}
       e^{-\gamma \tau}\ell(\xi(\tau))\geq \psi(x_0) + \frac{\varepsilon_2}{2} = V_{ \gamma}(x_0) + \frac{\varepsilon_2}{2},
   \end{equation}
    and there exists $\bar{\lambda}\in \Lambda$ such that
    \begin{align*}
       \varepsilon_3 &\leq \inf_{\csig\in \cfset} D_x\psi(x_0) \cdot f(x_0, \csig(0), \bar{\lambda}[\csig](0)) - \gamma \psi(x_0),\\
       &\leq D_x\psi(x_0) \cdot f(x_0, \csig(0), \bar{\lambda}[\csig](0)) - \gamma \psi(x_0),
    \end{align*}
    for all $\csig \in \cfset$. Then, there exists small enough $\delta_2 > 0$ such that
    \begin{align*}
       D_x\psi(\xi_{x_0}^{\csig, \bar{\lambda}}(\tau)) &\cdot f(\xi_{x_0}^{\csig, \bar{\lambda}}(\tau), \csig(\tau), \bar{\lambda}[\csig](\tau)) \\
       & - \gamma \psi(\xi_{x_0}^{\csig, \bar{\lambda}}(\tau))\geq  \frac{\varepsilon_3}{2}, %\forall \tau\in [0, \delta].
    \end{align*}
    for all $\tau\in [0, \delta_2]$.
    Similarly, we multiply both sides of the above inequality by $e^{-\gamma \tau}$ and integrate from $0$ to $\delta=\min\{\delta_1, \delta_2\}$ to get
    \begin{align*}
       e^{-\gamma \delta}\psi(\xi_{x_0}^{\csig, \bar{\lambda}}(\delta)) - \psi(x_0) \geq \frac{\varepsilon_3}{2}\delta.
    \end{align*}
    Recalling that $V_{\gamma}-\psi$ attains a local minimum of $0$ at $x_0$, we have
    \begin{align}\label{eq:sup_V>=eps}
       e^{-\gamma \delta}V_{\gamma}(\xi_{x_0}^{\csig, \bar{\lambda}}(\delta)) \geq  V_{\gamma}(x_0) + \frac{\varepsilon_3}{2}\delta.
    \end{align}
    Incorporating \eqref{eq:sup_l>=V} and \eqref{eq:sup_V>=eps} into the dynamic programming principle \eqref{eq: DPP}, we have 
    \begin{align*}
        V_{ \gamma}(x_0)&\geq \sup _{\lambda \in \Lambda}\inf_{\csig\in \cfset}
        \bigl\{\max[e^{-\gamma \delta} V_{\gamma}(\xi(\delta)), \min _{t\in [0, \delta]} e^{-\gamma t} \ell(\xi(t))]\bigr\},\\
        & \geq V_{\gamma}(x_0) + \max \left\{\frac{\varepsilon_2}{2}, \frac{\varepsilon_3}{2}\delta\right\},
    \end{align*}
    which is a contradiction, since $\varepsilon_2, \varepsilon_3, \delta > 0$. Therefore, we conclude that \eqref{eq: sup sol} must be true and hence $V_{ \gamma}$ is indeed a supersolution of \eqref{eq: HJI-VI}.
    
    For the uniqueness, we prove a comparison principle for our HJI-VI \eqref{eq: HJI-VI}, using similar techniques in the classical comparison and uniqueness theorems (see Theorem 2.12 in \cite{Bardi_Capuzzo-Dolcetta_2009}). We defer this proof to the appendix, and here show the uniqueness from another perspective: the contraction mapping.
    \end{proof}
\end{theorem}

% In fact, the uniqueness property can be seen as a result of the contraction property of the Bellman backup associated with the dynamic programming principle of \(V_{\gamma}\) in \eqref{eq: DPP}. 
These two theorems provide one numerical approach to compute the RA value function: solving the HJI-VI~\eqref{eq: HJI-VI}. This is not the only way of computing it. As will be shown in the next few paragraphs, the DPP induces a contractive Bellman backup operator. Therefore, value iteration can also be used to compute the RA value function.

We define a Bellman backup operator \(B_{T}\) : \(\operatorname{BUC}\left(\mathbb{R}^{n}\right) \mapsto \operatorname{BUC}\left(\mathbb{R}^{n}\right)\) for \(T>0\), (where \(\operatorname{BUC}\left(\mathbb{R}^{n}\right)\) represents the set of bounded and uniformly continuous functions: \(\mathbb{R}^{n} \mapsto \mathbb{R}\),) as
\begin{align} \label{eq: bellman_operator}
B_{T}[V](x) &\coloneqq \sup _{\lambda \in \Lambda}\inf_{\csig\in \cfset} \min \bigl\{ \notag\\
    &\min_{t\in [0, T]} \max [e^{-\gamma t}\ell(\xi(t)), \max_{s \in[0, t]} e^{-\gamma s} c(\xi(s))], \notag\\
    &\max[e^{-\gamma T} V(\xi(T)), \max _{t\in [0, T]} e^{-\gamma t} c(\xi(t))]\bigr\}.
\end{align}

With the help of the discount factor, we can show that this operator is a contraction mapping.
\begin{theorem}\label{th: contraction}
    (Contraction mapping). For any \(V^{1}, V^{2} \in \operatorname{BUC}\left(\mathbb{R}^{n}\right)\),
    \begin{equation}
            \left\|B_{T}\left[V^{1}\right]-B_{T}\left[V^{2}\right]\right\|_{L^{\infty}} \leq e^{-\gamma T}\left\|V^{1}-V^{2}\right\|_{L^{\infty}}, \label{eq: contraction}
    \end{equation}
    and the RA value function \(V_{\gamma}\) in \eqref{eq: def-inf} is the unique fixed-point solution to \(V_{\gamma}=B_{T}\left[V_{\gamma}\right]\) for each \(T>0\). Also, for any \(V \in \operatorname{BUC}\left(\mathbb{R}^{n}\right)\),
    \begin{equation}\label{eq: Bellman_limit}
        \lim _{T \rightarrow \infty} B_{T}[V]=V_{\gamma}.        
    \end{equation}
    
\begin{proof}
    Define
    \begin{equation*}
        w\left(\lambda, \csig, x\right)\coloneqq \max _{t\in [0, T]} e^{-\gamma t} c(\xi(t)), \ w^{i}\left(\lambda, \csig, x\right)\coloneqq e^{-\gamma T} V^{i}(\xi(T))
    \end{equation*}
 for \(i=1,2\). Then,
\begin{align*}
    B_{T}[V^{i}](x) &\coloneqq \sup _{\lambda \in \Lambda}\inf_{\csig\in \cfset} \min \bigl\{ \notag\\
    &\min _{t\in [0, T]} \max [e^{-\gamma t}\ell(\xi(t)), \max _{s \in[0, t]} e^{-\gamma s} c(\xi(s))], \notag\\
    &\max[w^{i}\left(\lambda, \csig, x\right), w\left(\lambda, \csig, x\right)]\bigr\}.
\end{align*}
Without loss of generality, %let $B_{T}\left[V^{1}\right](x) \geq B_{T}\left[V^{2}\right](x)$. 
let $V^{1}(x) \geq V^{2}(x)$. For any $\varepsilon>0$, there exists \(\bar{\lambda}\), such that
\begin{align*}
    B_{T}\left[V^{1}\right](x) - \varepsilon< &\inf_{\csig\in \cfset}  \min \bigl\{ \max[w^{1}\left(\bar{\lambda}, \csig, x\right), w\left(\bar{\lambda}, \csig, x\right)],\\
    \min _{t\in [0, T]} &\max [e^{-\gamma t}\ell(\xi(t)),
     \max _{s \in[0, t]} e^{-\gamma s} c(\xi(s))]\bigr\}, \\
    < &\min \bigl\{ \max[w^{1}\left(\bar{\lambda}, \csig, x\right), w\left(\bar{\lambda}, \csig, x\right)],\\
    \min _{t\in [0, T]} &\max [e^{-\gamma t}\ell(\xi(t)),
     \max _{s \in[0, t]} e^{-\gamma s} c(\xi(s))]\bigr\},
\end{align*}
for any $\csig$, which indicates that both of the following hold:
\begin{align}
    B_{T}\left[V^{1}\right](x) - \varepsilon <  &\max [w\left(\bar{\lambda}, \csig, x\right), w^{1}\left(\bar{\lambda}, \csig, x\right)], \label{eq: contraction BT V^1_1}\\
    B_{T}\left[V^{1}\right](x) - \varepsilon <  &\min _{t\in [0, T]} \max [e^{-\gamma t}\ell(\xi(t)),\notag \\
    &\quad \max _{s \in[0, t]} e^{-\gamma s} c(\xi(s))]. \label{eq: contraction BT V^1_2}
\end{align}

On the other hand, for any $\varepsilon>0$ and any $\lambda$, there exists $\bar{\csig}$ s.t.
\begin{align*}
    B_{T}\left[V^{2}\right](x) + \varepsilon > \min \bigl\{ \max[w^{2}\left(\lambda, \bar{\csig}, x\right), w\left(\lambda, \bar{\csig}, x\right)],\\
    \min _{t\in [0, T]} \max [e^{-\gamma t}\ell(\xi(t)),
     \max _{s \in[0, t]} e^{-\gamma s} c(\xi(s))]\bigr\},
\end{align*} 
 which indicates that either of the following holds:
\begin{align}
    B_{T}\left[V^{2}\right](x) + \varepsilon > &\max[w\left(\lambda, \bar{\csig}, x\right), w^{2}\left(\lambda, \bar{\csig}, x\right)], \label{eq: contraction BT V^2_1}\\
    B_{T}\left[V^{2}\right](x) + \varepsilon > &\min _{t\in [0, T]} \max [e^{-\gamma t}\ell(\xi(t)), \notag\\
    &\quad \max _{s \in[0, t]} e^{-\gamma s} c(\xi(s))]. \label{eq: contraction BT V^2_2}    
\end{align}
%Suppose first that $\max \left\{l\left(\bar{\lambda}, \bar{\csig}, x\right), l^{2}\left(\bar{\lambda}, \bar{\csig}, x\right)\right\} < \varepsilon + B_{T}\left[V^{2}\right](x)$.
Combining Eq.\eqref{eq: contraction BT V^1_1} and \eqref{eq: contraction BT V^2_1}, we have
\begin{align*}
B_{T}\left[V^{1}\right](x) &- B_{T}\left[V^{2}\right](x) \notag\\
&< 2 \varepsilon+\max \left\{w\left(\bar{\lambda}, \bar{\csig}, x\right) w^{1}\left(\bar{\lambda}, \bar{\csig}, x\right)\right\}\notag \\
&\quad -\max \left\{w\left(\bar{\lambda}, \bar{\csig}, x\right), w^{2}\left(\bar{\lambda}, \bar{\csig}, x\right)\right\}, \notag \\
&\leq  2 \varepsilon+\left|w^{1}\left(\bar{\lambda}, \bar{\csig}, x\right)-w^{2}\left(\bar{\lambda}, \bar{\csig}, x\right)\right|, \notag \\
&\leq  2 \varepsilon+e^{-\gamma T} \sup _{x \in \mathbb{R}^{n}}\left|V^{1}(x)-V^{2}(x)\right|.    
\end{align*}
The second inequality holds since, for all \(a, b, c \in \mathbb{R}, \mid \max \{a, b\}-\) \(\max \{a, c\}\left|\leq|b-c|\right.\). Moreover, combining Eq.\eqref{eq: contraction BT V^1_2} and \eqref{eq: contraction BT V^2_2}, we have $B_{T}\left[V^{1}\right](x) - B_{T}\left[V^{2}\right](x)< 2 \varepsilon$, so the result follows. As the above inequality holds for all \(x \in \mathbb{R}^{n}\) and \(\varepsilon>0\), we have
\[
\left\|B_{T}\left[V^{1}\right]-B_{T}\left[V^{2}\right]\right\|_{L^{\infty}\left(\mathbb{R}^{n}\right)} \leq e^{-\gamma T}\left\|V^{1}-V^{2}\right\|_{L^{\infty}\left(\mathbb{R}^{n}\right)}
.\]

Since \(V_{\gamma}\) is a fixed-point solution for all \(T>0\), the Banach's contraction mapping theorem \cite{evans2010chapter9.2} implies that \(V_{\gamma}\) is the unique fixed-point solution to \(B_{T}\left[V_{\gamma}\right](x)=V_{\gamma}(x)\) for all \(T>0\). In addition, we have
\[
\left\|B_{T}[V]-V_{\gamma}\right\|_{L^{\infty}\left(\mathbb{R}^{n}\right)} \leq e^{-\gamma T}\left\|V-V_{\gamma}\right\|_{L^{\infty}\left(\mathbb{R}^{n}\right)}.
\]
for all \(V \in \operatorname{BUC}\left(\mathbb{R}^{n}\right)\), thus we conclude \eqref{eq: Bellman_limit}.
\end{proof}
\end{theorem}

\begin{remark} \label{remark:uniqueness}
    Note that Theorem 3 and Theorem 4 both provide “uniqueness’’ results, but for two interconnected yet different objectives and therefore should not be confused. In Theorem 3, “uniqueness’’ is a property of the HJI-VI (16): we prove that there exists no other bounded Lipschitz continuous function that satisfies (16) in the viscosity sense, so the value function is the unique viscosity solution to (16). In contrast, Theorem 4 concerns the dynamic-programming side: it shows that the Bellman operator $B_{T}[V]$ in \eqref{eq: bellman_operator} is contractive, and hence, if the operator is applied recursively starting from an arbitrary initialization, the iteration converges to the unique fixed-point solution of $V = B_T[V]$. These two results give two independent numerical routes for computing the same value function: one via solving the HJI-VI backward in time, and one via value iteration induced by the contractive Bellman operator. They are related but conceptually distinct.
    %\zgrvs{It should be noted that Theorem~\ref{th: HJI} and Theorem~\ref{th: contraction} both provide ``uniqueness'' results. They study different but interconnected objectives and should not be confused. In Theorem~\ref{th: HJI}, ``uniqueness'' is the property of the HJI-VI~\eqref{eq: HJI-VI}. We prove that there exists no other bounded Lipschitz continuous function that satisfies~\eqref {eq: HJI-VI} in the viscosity sense. Theorem~\ref{th: contraction} shows that the Bellman operator~\eqref{eq: bellman_operator} is contractive; therefore, if the operator is applied recursively on a random initialization, it should converge to the same fixed point solution. These are distinct results.}
    % \zgrvs{On the other hand, the HJI-VI~\eqref{eq: HJI-VI} is merely the differential form of the DPP~\eqref{eq: DPP}, and the Bellman operator~\eqref{eq: bellman_operator} is derived based on the DPP~\eqref{eq: DPP}. Therefore, both result stems from the same equation, and should have inner connections. However, it is not clear to us whether those results are equivalent, that is, the HJI-VI has a unique viscosity solution if and only if the corresponding Bellman operator is contractive.}
\end{remark}
\subsection[Computation of V gamma]{Computation of $V_\gamma$}
Theorem \ref{th: contraction} allows the use of various methods to compute the value function $V_{\gamma}$ using the operation $B_{T}[\cdot]$, which does not require a specific initialization. For instance, a popular method is to solve the finite-horizon HJ equation presented in the following proposition. The proof is analogous to the proof of Lemma 5 of~\cite{choi2023forward}.

\begin{proposition}\label{prop:unique_VS_finite}(Finite horizon HJI-VI for the computation of $V_{\gamma}$). For a given initial value function candidate \(V^{0} \in \operatorname{BUC}\left(\mathbb{R}^{n}\right)\), and let \(W\) : \(\mathbb{R}^{n} \times  [0, T]  \mapsto \mathbb{R}\) be the unique viscosity solution to the following terminal-value HJI-VI
   \begin{align}
        &W(x, T) = \max\{\ell(x), c(x), V^0(x)\}, \forall x\in \R^n, \label{eq: finite_init}\\
        &0=\max\bigl\{\min \{D_tW(x, t) + \max _{d \in \dset}\min_{u\in \cset} D_xW(x, t) \cdot f(x, u, d)  \notag\\
        & \qquad - \gamma W(x, t),\ell(x)-W(x, t)\}, c(x)-W(x, t)\bigr\},  \label{eq: finite HJI-VI}
   \end{align}   
for \((x, t)\in \R^{n}  \times [0, T]\). 

Then we have: \(W(x, 0) \equiv B_{T}\left[V^{0}\right](x)\).
\end{proposition}
In Prop. \ref{prop:unique_VS_finite}, any $V^0\in\textnormal{BUC}(\R^n)$ works for the computation of $V_\gamma$; for instance, a straightforward choice of $V^0$ can be $c(x)$. As $T\rightarrow\infty$, $D_tW(x, t)$ vanishes to 0 for all $x\in\R^n$. 

Prop. \ref{prop:unique_VS_finite} emphasizes that given any initialization $V^0$, the results from the contraction mapping (Theorem~\ref{th: contraction}) and from solving the terminal value HJI-VI are the same, if the terminal value of the HJI-VI satisfies ~\eqref{eq: finite_init}.  
% \begin{remark} Although any $V^0\in\textnormal{BUC}(\R^n)$ works for Lemma \ref{lemma:numerical_finite_HJPDE} and for the computation of $\VV$, the initial value of the solution to \eqref{eq:numerical_finite_HJPDE} has to be greater than or equal to $h_S$ to satisfy \eqref{eq:numerical_finite_HJPDE} at $t=0$. This is why the initial value is set as \eqref{eq:numerical_finite_HJPDE_initialCondition}. A straightforward choice of $V^0$ can be $h_S$. Also, as $T\rightarrow\infty$, $\frac{\partial W}{\partial t}$ vanishes to 0 for all $x\in\R^n$.
% \end{remark}

Combining Theorem \ref{th: contraction} and Prop. \ref{prop:unique_VS_finite}, we have 
\begin{align}
    \lim_{T\rightarrow \infty}B_T[V^0]=\lim_{T\rightarrow \infty} W(x,0)=V_\gamma(x).
\end{align}
The PDE \eqref{eq: finite HJI-VI} can be numerically solved backward in time from the terminal condition \eqref{eq: finite_init}, by using well-established time-dependent level-set methods \cite{Mitchell2005b}.  

%Theorem 4 facilitates the implementation of various numerical methods based on time discretization, such as value iteration, to approximate the solution of \( V_{\gamma} \). The subsequent corollary of Theorem 4 establishes that value iteration, initialized with any function \( V^0 \in \text{BUC}(\mathbb{R}^n) \), will converge to \( V_{\gamma} \) with a Q-linear convergence rate, as specified in Equation (33). For a given time step \( \Delta t \), the semi-Lagrangian approach can be utilized to approximate the exact Bellman operator in Equation (27), leading to a numerical approximation. Furthermore, as \( \Delta t \to 0 \), the resulting value function converges to \( V_{\gamma} \).
In addition, another line of methods enabled by Theorem \ref{th: contraction} is based on time discretization, such as value iteration, to accurately solve for \( V_{\gamma} \). The subsequent corollary of Theorem \ref{th: contraction} establishes that value iteration, initialized with any function \( V^0 \in \text{BUC}(\mathbb{R}^n) \), will converge to \( V_{\gamma} \) with a Q-linear convergence rate, as specified in \eqref{eq: value iteration}. For a given time step \( \Delta t \), the semi-Lagrangian approach can be utilized to approximate the exact Bellman operator in \eqref{eq: contraction}, leading to a numerical approximation. Furthermore, as \( \Delta t \to 0 \), the resulting value function converges to \( V_{\gamma} \) \cite{Kene_2018}.
\begin{corollary}(Value Iteration). \label{corollary:value_iteration} Let \( V^0 \in \text{BUC}(\mathbb{R}^n) \) and consider a time step \( \Delta t > 0 \). Define the sequence \( \{V^k\}_{k=0}^{\infty} \) iteratively as $V^k := B_{\Delta t}[V^{k-1}], \text{ for } k \in \mathbb{N}$. Then, the following holds:
\begin{align}
\frac{\|V^{k+1} - V_{\gamma}\|_{\infty}}{\|V^k - V_{\gamma}\|_{\infty}} \leq e^{-\gamma \Delta t} < 1,\label{eq: value iteration}
\end{align}
which implies that $\lim_{k \to \infty} V^k = V_{\gamma}$.

\begin{proof}
This result follows directly from Theorem \ref{th: contraction}.
\end{proof}
\end{corollary}

\subsection{The Finite-horizon Value Function and Control Synthesis}
\label{sec: controller syns}
In many HJ-reachability based works, the optimal controller can be synthesized by the gradient of the value function (see Sec I.5 in~\cite{Bardi_Capuzzo-Dolcetta_2009}). However, it should be pointed out that the controller synthesized directly from the gradient of ~\eqref{eq: def-inf} does not guarantee finite-time reach-avoid of the target and constraint sets. Though from Prop~\ref{prop: def-RA}, we concluded that $V_\gamma(x)< 0$ implies there exists some $T$ s.t. $\ell(\xi(T))<0$ and $c(\xi(t))< 0 $ for all $t \in [0,T]$, it does not mean the optimal control determined from the gradient of $V_\gamma(x)$ is able to drive the system to the target at exactly time $T$. One reason is that HJI-VI~\eqref{eq: HJI-VI} does not provide any useful information ($\max_d \min_u \dot V_\gamma = \gamma V_\gamma$ has nothing to do with finite-time reach avoid), and that the system can travel freely before $T$, and then apply some control and reach the target. 

In this section we discuss how to design a time-optimal RA controller using a finite-horizon version of the infinite RA value function discussed above. %Similarly to its  counterpart, W(x,t)$ is defined below.

\begin{definition}
A finite-horizon RA value function $W(x, t): \mathbb{R}^{n} \times  [0, T]  \rightarrow \R$ is defined as
    \begin{align}\label{eq: def-finite}
        W(x, t) &\coloneqq \sup _{\lambda \in \Lambda}\inf_{\csig\in \cfset}\min_{\tau \in[t, T]} \notag\\
        & \hspace{-2em}\max\bigl\{e^{-\gamma (\tau - t)} \ell(\xi^{\csig, \dmap[\csig]}_{x, t}(\tau)), \max_{s\in [t, \tau]}e^{-\gamma (s-t)} c(\xi^{\csig, \dmap[\csig]}_{x, t}(s))\bigr\}.
    \end{align}    
\end{definition}

The zero sublevel set of $W(x,t)$ characterizes a finite-horizon RA set, i.e., it is the set of initial states that can reach the target at some $\tau \in [0,t]$, and avoid the obstacle in $[0,\tau]$. We state without proof that~\eqref{eq: def-finite} satisfies the corresponding DPP and is the unique viscosity solution to the HJI-VI~\eqref{eq: finite_init}~\eqref{eq: finite HJI-VI}, by taking $ V^0(x) =\max\{\ell(x), c(x)\} $.

% It satisfies the corresponding DPP and is the unique viscosity solution to \eqref{eq: finite HJI-VI}. 

From the definition, for any fixed time horizon $T$, $ W(x, 0) \geq V_{\gamma}(x) $. This means the zero sub-level set of $ W(x, 0)$ provides an under-approximation of the infinite-time RA set. If we take $T \rightarrow \infty$, the finite-time RA value function at $t = 0$ is exactly the same as~\eqref{eq: def-inf}, i.e., $ W(x, 0) = V_{\gamma}(x) $. Further, for a fixed state $x$ and time horizon $T$, as $t$ decreases from $T$ to 0, $W(x,t)$ is non-increasing, and the time (if exists) $W(x,t) $ first decay to $0$ is the minimal time for the trajectory starting from $x$ to reach the target while avoiding the obstacle. 
% Denote that time as $t^*\coloneqq \min\{t: V_\gamma(x, T-t) = 0\}$, 
One optimal control signal along a trajectory $\xi(t)$ is given by
% We defer the discussion to appendix and here focus on controller synthesis. 
    % The optimal controller is then given by
    \begin{equation}\label{eqn:RA_ctrl}
        \pi_{RA}(\xi(t))=\argmax\limits_{\dstb \in \dset}\min\limits_{\csig \in \cset}D_x W(\xi(t), T- t^*) \cdot f(\xi(t), u, d),
    \end{equation}
where $t^*\coloneqq \min\{t: W(\xi(t), T-t) = 0\}$. Notice here $\xi(t)$ and $\pi_{RA}(\xi(t))$ is an optimal control-trajectory pair~\cite{Bardi_Capuzzo-Dolcetta_2009}. 

%It turns out that $V_\gamma(x,t)$ is the unique viscosity solution to \eqref{eq: finite HJI-VI}. Let $t^*\coloneqq \min\{t: V_\gamma(x, T-t) = 0\}$.
%In the previous section we introduced an infinite-horizon RA value function whose zero sublevel set characterizes the $\mathcal{R} \mathcal{A}(\mathcal{T}, \mathcal{C})$

%\begin{remark}
%    It should be noted that in Proposition~\ref{prop:unique_VS_finite}, the HJ-PDE has a unique solution. The result can be derived from the standard comparison principle. However, in Theorem~\ref{th: HJI}, the solution to the PDE is not unique. Though in~\cite{}, the authors claim that with a discount factor added to the value function and the PDEs, the solution to the PDE become unique, their PDE indeed has non-unique solutions. It is also the case for~\eqref{eq: def-inf} For counter examples, check~\cite{}. \zgnote{cite Kene's paper for the discount factor thing, put the counter examples on arxiv and cite it.} 
%\end{remark}

\subsection{Impact of the Discount Factor \texorpdfstring{$\gamma$}{gamma}}\label{subsec:impact of gamma}

Theoretically, there is no harm choosing $\gamma$ as large as possible. Discounting guarantees boundedness, uniqueness of solutions to the HJI-VI \eqref{eq: HJI-VI}, and $V_\gamma$ being globally Lipschitz if the dynamics satisfies $L_f<\gamma$. Remark~\ref{remark:same_zero_set} shows that for any $\gamma>0$ we recover the same certified set, $ \mathcal{R} \mathcal{A}(\mathcal{T}, \mathcal{C})=\{x:\,V_\gamma(x)<0\}$, so changing $\gamma$ does not alter the RA set itself. 

Nevertheless, since we rely on a numerical solver to solve the value functions, it is necessary to consider the impact of $\gamma$ from a numerical perspective. The numerical scheme we use to compute the RA value function is analogous to~\cite{Fisac_Chen_Tomlin_Sastry_2015}. We discretize the state space into grids and the time into a discrete time series. Two direct results from this discretization are 1) the state-space and time series will be both finite and bounded, and 2) the signed-distance function can always be used as the cost and constraint functions (which are bounded and Lipschitz continuous). For each grid point, we compute the finite-time RA value~\eqref{eq: def-finite} from $t=0$ to a prespecified time horizon $T$. We further set a convergence threshold and keep track of the state-wise change of value at each timestep. If the change of value is less than the threshold for all grid points at some time $t$, we say the value function has converged, and use $W(x,t)$ to approximate $V_\gamma(x)$. 

With this approximation scheme, the discount factor $\gamma$ plays a numerical—rather than semantic—role for the reach--avoid certification. That is, changing $\gamma$ affects the numerical values of $V_\gamma$ and the time to compute them, but not the certified sets. In particular, for Value Iteration, a larger $\gamma$ strengthens the contraction of the Bellman backup operator and thus accelerates the computation. This is the benefit of choosing a moderately large $\gamma$. 

On the other hand, setting $\gamma$ excessively large can introduce practical numerical issues. In our discrete scheme, for the states whose values are initialized as negative, since~\eqref{eq: def-inf} keeps track of the infimum over $t$, the values will not increase as computation proceeds. 
% Therefore, the initial values become an upper bound for the intermediate values of this state. 
For grid points that start with a positive value but actually belong to the RA set, a very large $\gamma$ pushes the value down to (almost) zero within only a few iterations. After that, the iteration-to-iteration change is at machine precision, so even with a tight stopping tolerance (for example, $O(10^{-8})$) the solver may stop too early, i.e., the value function may ``falsely converge.'' As will be shown in the 1D example, the value function ``converges'' within a few iterations with a very large $\gamma$, but it does not recover the true RA set. Without further information, we cannot tell if the value function truly converges. %For the states initialized with a positive value, $\gamma$ will decrease the value to zero in the first few iterations (if the states are in the RA set), and afterward, the change of value at each time step will be very tiny. This means even if we set a very small} \zgrvs{convergence threshold (e.g., $O(10^{-8})$), the value function may ``falsely converge.'' As will be shown in the 1D example, the value function ``converges'' within a few iterations with a very large $\gamma$, but it does not recover the true RA set. Without further information, we have no idea if the value function truly converges. }

This ``false convergence'' is problematic for controller synthesis. Suppose an initial state needs 10 timesteps to reach the target, then from~\eqref{eqn:RA_ctrl}, we need at least the intermediate value functions $W(x,t)$ from timesteps 1 to 10. However, with a large $\gamma$, the computation may ``falsely converge'' within 5 timesteps. If this happens, we do not have enough time snapshots of the value function to derive the control.%i.e., we need $W(x,t)$, where $t= 1,2,...,10$. But with a large $\gamma$, it may ``falsely converge'' within 5 timesteps, i.e., the computation will stop at $W(x,5)$. If this happens, we do not have enough time snapshots of the value function to derive the control.

Besides ``false convergence", an overly large $\gamma$ also amplifies numerical errors near the boundary of the RA set, and we call this phenomenon ``boundary blow-up''. In a grid-based HJI-VI solver, the value at one node is influenced by the neighboring nodes. When the true value inside the RA/SA set is a very small negative number (say $O(10^{-8})$), the neighboring states just outside the set can still have noticeably positive values (say $O(10^{-4})$). A large $\gamma$ makes this contrast sharper, and the positive values can dominate the update. As a result, boundary points that should be negative are computed as positive, and the zero level set becomes inaccurate.

%\zgrvs{Besides the ``false convergence'', a very large $\gamma$ can cause another problem: the numerical value function may have relatively larger errors near the boundary of the RA set. This is because numerically, to compute the value of $x$, the solver relies on the value of nearby states of $x$. With a very large $\gamma$, if we focus on the states in the RA set that are near the boundary, we know in theory the value should be negative, but it might be a very tiny negative value (e.g., $-10^{-8}$). However, the numerically solved value might be dominated by the positive values (e.g., $10^{-4}$) outside the RA set, and therefore becomes positive.  }

\begin{figure}%[t]
    \centering
   % \vspace{.5em}
    \includegraphics[width=0.49\textwidth]{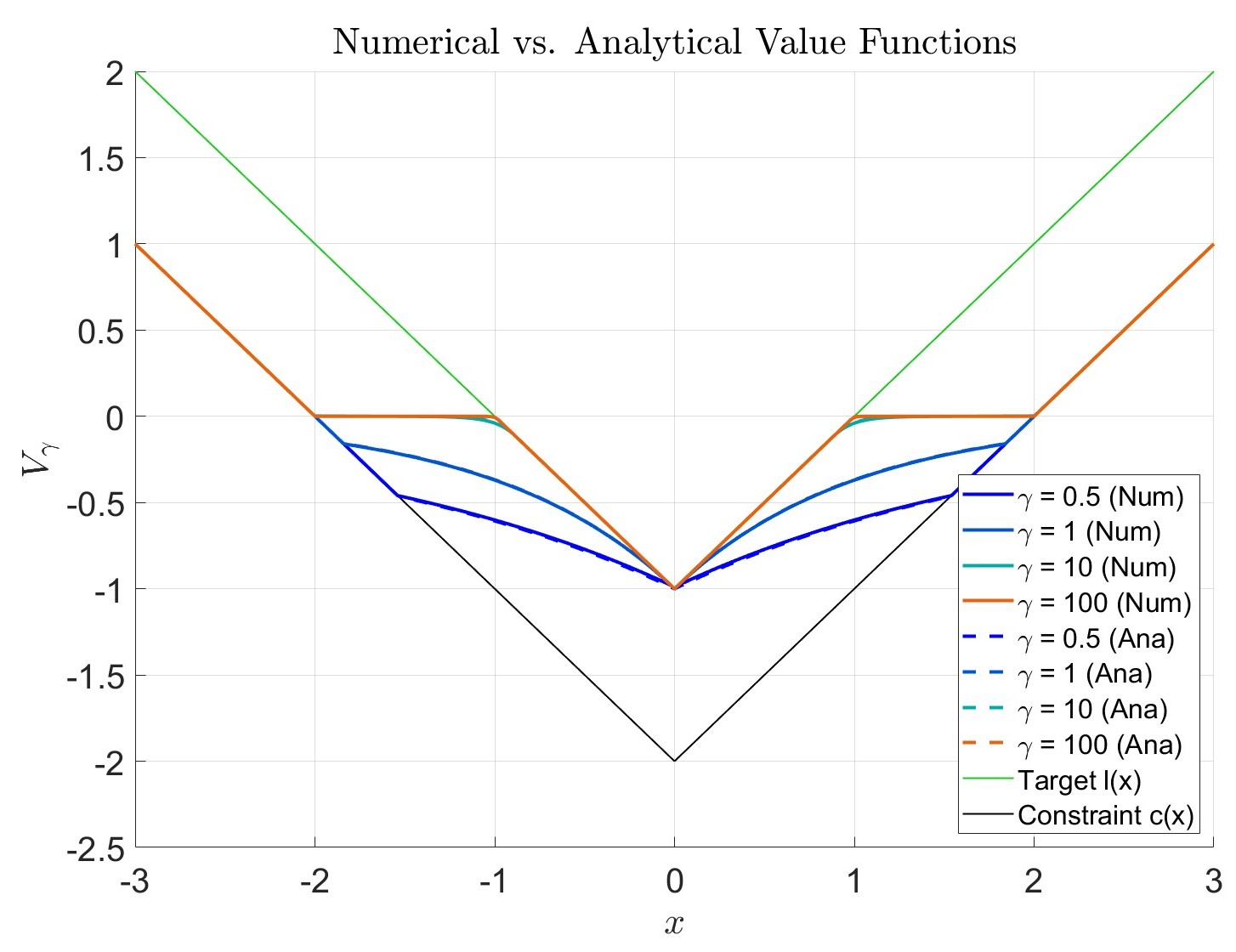}
    \caption{1D RA value function $V_\gamma(x)$ and its level sets. Solid lines are the numerical solutions, and the dashed lines are the analytic solutions. For all different $\gamma>0$, the zero-level sets are the same. For all $\gamma$'s, the analytical solution and the numerical solution are almost identical. The computation time (averaged for 20 runs) for $\gamma=0.5,1,10, 100$ is 78 timesteps, 85 timesteps, 28 timesteps, and 6 timesteps, with convergence threshold $5\times10^{-4}$. Only $\gamma=0.5,1$ provides enough snapshots of the value function to derive the control. When $\gamma = 100$, the change of value decreases to 0 (under MATLAB precision) after 14 iterations. }
 %   \vspace{-1em}
    \label{fig:V_1D_RA}
 %   \vspace{-1em}
\end{figure}

We now use a simple 1-dimensional running example to show the impact of the discount factor $\gamma$. We will derive the analytical form of the RA value function, where $\gamma$ is a hyperparameter, and compare with the numerical solutions. Let's consider the following system:
\begin{align} \label{eq:ex_sys}
    \dot x = u, \hspace{1em} u \in [-1,1].
\end{align}
Assume we want the system to reach the target $\target = \{ x\colon x\in(-1,1) \}$, and the constraint $\obs = \{ x\colon x\in(-2,2) \} $, so we set $\ell(x) = |x| - 1$, $c(x) = |x|-2$. Under this setup, the optimal feedback control is given by:
\begin{align}\label{eq:ex_ctrl}
    \pi = \begin{cases}
        -1 & x>0 \\1& x<0 \\ 0 &x=0
    \end{cases}.
\end{align} For this system, the RA set is given by the interval $(-2,2)$, and it takes $1s$ for initial states $x = 2^-$ and $-2^+$ that are very close to the boudaries to reach the target. If we fix the time step to be $0.02s$, this means we need at least 50 snapshots of the value function (i.e., $W(x,t)$, where $t= 1,2,...,50$) to generate the trajectory. Notice that for this system, the Lipschitz constant is $L_f = 0$, therefore, $\gamma$ can be any positive value. The RA value function for this problem depends on the value of $\gamma$. More specifically, we have
\begin{align} \label{eq:ex_RA_val}
     &V_\gamma(x) \notag \\
     =& \begin{cases}
         \max \bigl \{  -e^{-\gamma |x|} , |x|-2   \bigr \} & \gamma\leq 1,\\
         \max \bigl \{ \min (|x|-1, -e^{-\gamma |x|} ), \\\hspace{7em}|x|-2   \bigr \} & \gamma >1, |x|<1-\frac{1}{\gamma}, \\
          \max \bigl \{ \min ( -\frac{1}{\gamma} e^{-\gamma |x| -1 + \gamma}, \\
          \hspace{3em}-e^{-\gamma |x|} ), |x|-2   \bigr \} & \gamma >1,  |x|\geq 1-\frac{1}{\gamma}.
     \end{cases} 
\end{align} The detailed derivation is presented in the appendix. Taking the derivative of~\eqref{eq:ex_RA_val}, it can be seen that all of the Lipschitz constants of the RA value functions for different $\gamma$ are smaller than 1 (this can be verified by taking the derivative of all components inside the $\min$ and $\max$ sign). This is proved by Proposition~\ref{prop: V_RA properties}: the Lipschitz constants of the RA value function depend on the Lipschitz constant of $\ell(x)$ and $c(x)$. Further, it can be verified that the zero sub-level sets for different $\gamma$'s are all $x \in (-2,2)$, verifying Remark~\ref{remark:same_zero_set}. Fig. \ref{fig:V_1D_RA} illustrates what we discussed above.

% \zgrvs{However, $\gamma$ does change the computation time of the value function. Empirically, a larger $\gamma$ provides a faster convergence. The comparison of the numerical solution and analytical solution~\eqref{eq:ex_RA_val} is shown in Fig. \ref{fig:V_1D_RA}.}

\begin{figure}%[t]
    \centering
   % \vspace{.5em}
    \includegraphics[width=0.49\textwidth]{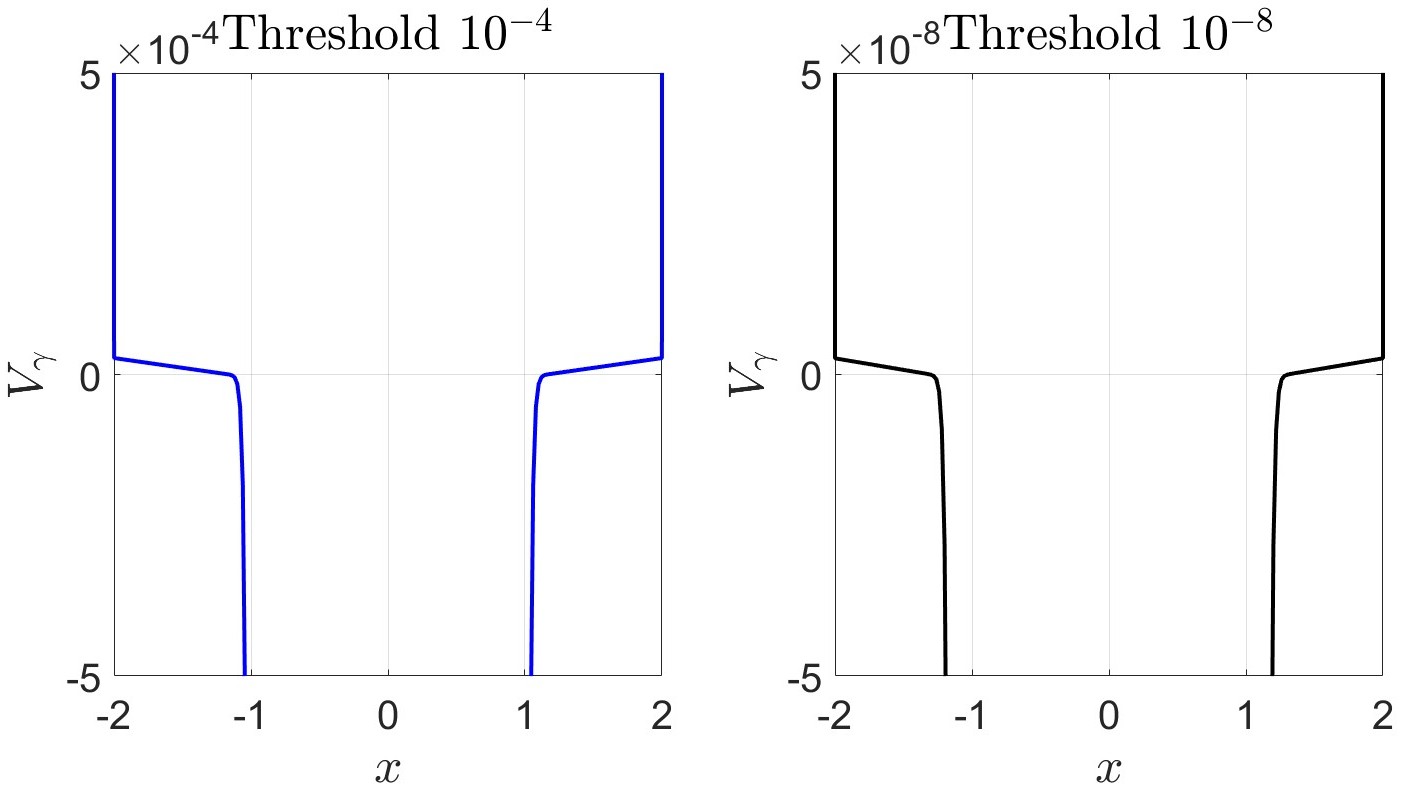}
    \caption{1D RA value function $V_\gamma(x)$ for $\gamma = 100$ with convergence threshold $10^{-4}$ (left) and $10^{-8}$ (right). For both cases, the values near $x = \pm 2$ are very small positive values, whereas in theory they should all be negative. This is because the nearby states have a relatively large positive value, and dominate the value inside the RA set. }
 %   \vspace{-1em}
    \label{fig:V_1D_RA_closeup}
 %   \vspace{-1em}
\end{figure}

Note that for $\gamma=100$, when $|x|\geq 2$, $V_\gamma(x) = |x|-2$. When $|x|<2$, $V_\gamma(x) =\min ( -\frac{1}{\gamma} e^{-\gamma |x| -1 + \gamma}, -e^{-\gamma |x|} ) <0 $. However, as can be seen in Fig.~\ref{fig:V_1D_RA_closeup}, the numerically solved value function remains positive (though the magnitude is tiny) near $x = \pm2$. This is because the relatively larger positive values outside the RA set dominate the value of those states.

\begin{remark}(Choice of the discount factor $\gamma$)\label{remark: effect of gamma}
In practice, $\gamma$ must be chosen carefully. It should be large enough to satisfy $L_f<\gamma$ and obtain a reasonable contraction rate, but small enough to avoid causing the issues explained above. A convenient guideline is to choose parameters such that $\gamma\in [L_f, 3L_f]$ if $L_f>0$ and $\gamma\in [0, 1]$ if $L_f=0$. We acknowledge that the study of numerical properties for accurately solving HJI-VI is beyond the scope of this work and will be left for future research.
\end{remark}

\section{TWO-STEP FORMULATION FOR STABILIZE-AVOID PROBLEM}\label{sec: SA}
In this section, we propose a two-step method to solve the SA problem by combining the RA formulation \eqref{eq: def-inf} and the R-CLVF. We further assume $\ell(x) $ is upper bounded by $r(x;p) = \| x- p\| - a$, i.e. $\ell(x) \leq r(x;p)$, and the zero sub-level set of $r(x)$ contains some robust control invariant subset $\mathcal I$. With this assumption, there exists at least one sub-level set of the R-CLVF to be a strict subset of the target set. Note also that this assumption is less strict than assuming there exists a stabilizable equilibrium point inside the target set.
% The constants $p$ and $a$ are defined in the following Definition \ref{def: R_CLVF}. %To be more specific, we first find the 
%Now we are ready to charaterize the SA set. 
%The main idea is to treat the largest robust control invariant set $I_\text{h}$ of $\clvf$ contained in the target $\target$ to be a new target set of the RA value function \eqref{eq: def-inf}. Let $h$ be the associated level of $I_\text{h}$, and define the shifted R-CLVF $\sclvf(x)\coloneqq \clvf(x) - h$. We construct a new discounted value function:

The idea is straightforward: we treat one level set of the R-CLVF as the new target set. There are two issues to be settled. The first is, which level set of the R-CLVF should we choose? Ideally, this level set should be the largest sub-level set $\mathcal I_\text{M}$ of the $\clvf$ contained in the target $\target$. 

The second issue relates to the properties of the R-CLVF itself. As shown in~\cite{gong2024robustcontrollyapunovvaluefunctions}, the R-CLVF is not globally bounded on $\mathbb R^n$. Therefore, if we simply replace $\ell$ with $\clvf$, the corresponding value function will not enjoy the same properties as the RA value function. In fact, any Control Lyapunov Function (CLF) is radially unbounded, and naively replacing $\ell$ with a CLF will have the same issue. The solution is straightforward; we could truncate the R-CLVF (and any valid CLF). Define the shifted R-CLVF:
\begin{align*}
    \sclvf (x) \coloneqq \begin{cases}
        \clvf (x) - M &  \clvf (x) \leq  N, \\
        N-M  &  \clvf (x) >  N.
    \end{cases}
\end{align*}
% To fit in our RA framework, we define the shifted R-CLVF $\sclvf \coloneqq \clvf - M$, 
where $M$ is the level of $\mathcal I_\text{M}$, $N$ is a user specified truncation parameter, and $N>M$. The shifted R-CLVF is guaranteed to be Lipschitz continuous and bounded on $\mathbb R^n$. Now we construct the SA value function $\VSA(x) \colon \mathbb R^n \mapsto \mathbb R$:
\begin{align}\label{eq: def-V^SA}
    \VSA(x) &\coloneqq \sup _{\lambda \in \Lambda}\inf_{\csig\in \cfset}\inf_{t\in[0, \infty)} \notag\\
    &\max\{e^{-\gamma t} \sclvf(\xi(t)), \sup_{s\in [0, t]}e^{-\gamma s} c(\xi(s))\},
\end{align}
where the cost function $\ell$ in \eqref{eq: def-inf} is replaced by $\sclvf$. The SA value function shares all the properties of the RA value function. In fact, it is indeed a RA value function whose target set is $\mathcal I_\text{M}$, where $\mathcal I_\text{M} = \{x\colon \sclvf(x) < 0 \} $. Though simply truncating a function may cause problems in some cases, it does not for us. This is because we truncate the R-CLVF on its positive part; therefore, it will only affect the positive part of the SA value function. Indeed, as will be shown in the next Proposition, the negative part of the SA value function is what we need.

Next, we show that its zero sublevel set is the desired SA set.
\begin{proposition}(Exact Recovery of SA set) \label{prop:exactSA}
% If we define
% \begin{align*}
%     \mathcal{S} \mathcal{A}(\mathcal{T}, \mathcal{C})\coloneqq& \bigl\{x \in \mathbb{R}^{n}\colon \forall \lambda \in \Lambda, \exists \csig\in \cfset \text { and } T \geq 0,\\  \text { s.t., } & \forall t\geq 0,
%     \xi_{x}^{\csig, \lambda[\csig]}(t) \in \constraint \text { and },\\ & \lim_{t\rightarrow \infty }\min_{y\in \partial \mathcal I_\text{m}} \| \xi_{x}^{\csig,  \lambda[\csig]}(t) -y \| =0 \bigr\},
% \end{align*}then the recovery is exact if $\partial\mathcal I_\text{m}\cap \mathcal I_\text{M} = \emptyset$.% i.e., $\mathcal I_\text{m}$ is a strict subset of $\mathcal I_\text{M}$.
    \begin{equation}\label{eq: def-SA}
        \mathcal{S} \mathcal{A}(\mathcal{T}, \mathcal{C}) = \{x\colon \VSA(x) < 0 \}    
    \end{equation}
%    \bynote{Find the proper definition of $\mathcal{S} \mathcal{A}(\mathcal{T}, \mathcal{C})$.\\}
\begin{proof}
     For sufficiency, note that $x\in  \mathcal{S} \mathcal{A}(\mathcal{T}, \mathcal{C})$ implies there exists $\csig$ such that for any $\dmap$, it safely steers the system to $\mathcal T$ and ultimately converge to $\mrcis$. Since the trajectory is continuous, it must enter $\mathcal I_\text{M}$, at some time $T\geq 0$, i.e., $\sclvf(\xi^{\csig, \dmap}_{x}(T))< 0$ and $c(\xi(s))<0$ for all $s\in[0,T]$. Since we take the infimum over $[0,\infty)$, and at $t = T$, the value is already negative, we conclude that $\VSA(x) < 0$. 
     
     % that there exists $T\geq 0$ such that $\xi(T)\in \mathcal I_\text{M}$, i.e., $\VSA < 0$. This proves the sufficiency. 
     % since  $\partial\mathcal I_\text{m}\cap \mathcal I_\text{M} = \emptyset$. Thus, by definition, $\VSA(x)<0$.
     
    Now suppose $\VSA(x)< 0$. Then, there exist some $T$, for all $\dmap_{[0,T]}$, there exists $\csig _{[0,T]}$  s.t 
    \begin{align*}
        \max \{e^{-\gamma T} \sclvf(\xi(T)), \sup_{s\in [0, T]}e^{-\gamma s} c(\xi(s))\}<0.
    \end{align*}
    Since exponential is positive, this means $ \sclvf(\xi(T)) < 0$ and $ c(\xi(s))<0$ for all $s \in [0,T]$. Further, as a direct result of Theorem.~\ref{thrm: clvf_exp_decay}, since $ \sclvf(\xi(T)) < 0$, for all $\dmap_{[T,\infty)}$, there exists $\csig _{[T,\infty)}$ s.t. $\lim_{t\rightarrow \infty }\min_{y\in \partial \mathcal I_\text{m}} \| \xi_{x}^{\csig,  \lambda[\csig]}(t) -y \| =0 $. This means for all $\dmap$, there exists $\csig$ that steers the systems to $\mathcal T$ and then
    % Next, applying $\pi_H$ stabilizes the system to $\mathcal{I}_m$. To conclude, the switching controller
    % \begin{equation}
    %     \pi_{\mathcal{SA}}(x)= \begin{cases}\pi_{RA}(x) & \text { if } \sclvf > 0 \\ \pi_{H}(x) & \text { otherwise }\end{cases}
    % \end{equation}
    safely stablize the system to $\mathcal{I_{\text{m}}}$. 
    % which means $x\in  \mathcal{S} \mathcal{A}(\mathcal{T}, \mathcal{C})$.
    %Then, there exists $\pi_{RA}$ that for any $\dsig$, it steers the system to the zero sublevel set of $\sclvf$ at time $t_1$, i.e., $\sclvf(\xi$
\end{proof}
\end{proposition}

%\zgrvs{If we make a stricter assumption that there exists a stabilizable equilibrium point $p_e$ inside the target set, then there should exist a control Lyapunov function (CLF) that stabilizes the system to the equilibrium point. Furthermore, picking $p = p_e$, the R-CLVF is a Lipschitz continuous CLF with respect to the equilibrium point $p_e$. One may obtain a CLF from other approaches: for simple low-dimensional systems, analytic CLFs can be obtained, neural CLF~\cite{NEURIPS2019_2647c1db,10015199} is another available option. Denoting the CLF as $V_{\text{CLF}}$ and replacing $\sclvf$ with $V_{\text{CLF}}$ in~\eqref{eq: def-SA}, one can still prove Proposition~\ref{prop:exactSA}. However, by making this assumption, the SA problem becomes restricted. Because the existence of a stabilizable equilibrium point inside the target set $\target$ is not guaranteed.}

%\begin{remark}
%    \zgrvs{The key to our formulation is to reach the region of stabilizability (also called the region of null-controllability), and then stabilize to some set while avoiding the obstacle. The R-CLVF serves as a perfect tool for our formulation, as it relaxes the assumption on the system dynamics and allows for arbitrarily picking the target set. }
%\end{remark}
As mentioned before, the SA problem is solved in a two-step manner: we first reach one level set of the R-CLVF ($\mathcal I_\text{M}$), then use the R-CLVF to stabilize the system. Therefore, one feedback controller can also be synthesized in the two-step formulation: given any initial states in the SA set, use $\pi_{RA}$~\eqref{eqn:RA_ctrl} in the reach-avoid phase, and $\pi_{H}(x)$~\eqref{eqn:opt_ctrl} in the stabilize-avoid phase:
\begin{equation} \label{eqn:SA_ctrl}
    \pi_{\mathcal{SA}}(x)= \begin{cases}\pi_{RA}(x) & \text { if } \sclvf > 0, \\ \pi_{H}(x) & \text { otherwise .}\end{cases}
\end{equation}

% \zgrvs{On the other hand, with R-CLVF, }
A key feature of our two-step SA framework is its modularity concerning the stabilization component. While this paper employs the R-CLVF for its broad applicability, the framework can readily incorporate any valid CLF. Specifically, under the stricter assumption that the target set $\target$ contains a known, stabilizable equilibrium point, $p_e$, a user-provided CLF that ensures convergence to $p_e$ can be substituted for the R-CLVF. Such CLFs may be derived from various approaches, including analytical methods for low-dimensional systems or learning-based techniques such as Neural CLFs~\cite{NEURIPS2019_2647c1db,10015199}. By denoting this function as $V_{\text{CLF}}$ and replacing the shifted R-CLVF ($\sclvf$) in the definition of the SA value function~\eqref{eq: def-V^SA}, the theoretical guarantees of our framework, including the exact recovery of the SA set stated in Proposition~\ref{prop:exactSA}, remain valid. This flexibility is noteworthy; however, the reliance on a pre-existing equilibrium point restricts the problem's scope, as such a point is not guaranteed to exist in general. Our use of the R-CLVF bypasses this limitation, offering a more broadly applicable solution.

\section{NUMERICAL EXAMPLES} \label{Sec:numerical}
\begin{figure*}%[t]
    \centering
   % \vspace{.5em}
    \includegraphics[width=\textwidth]{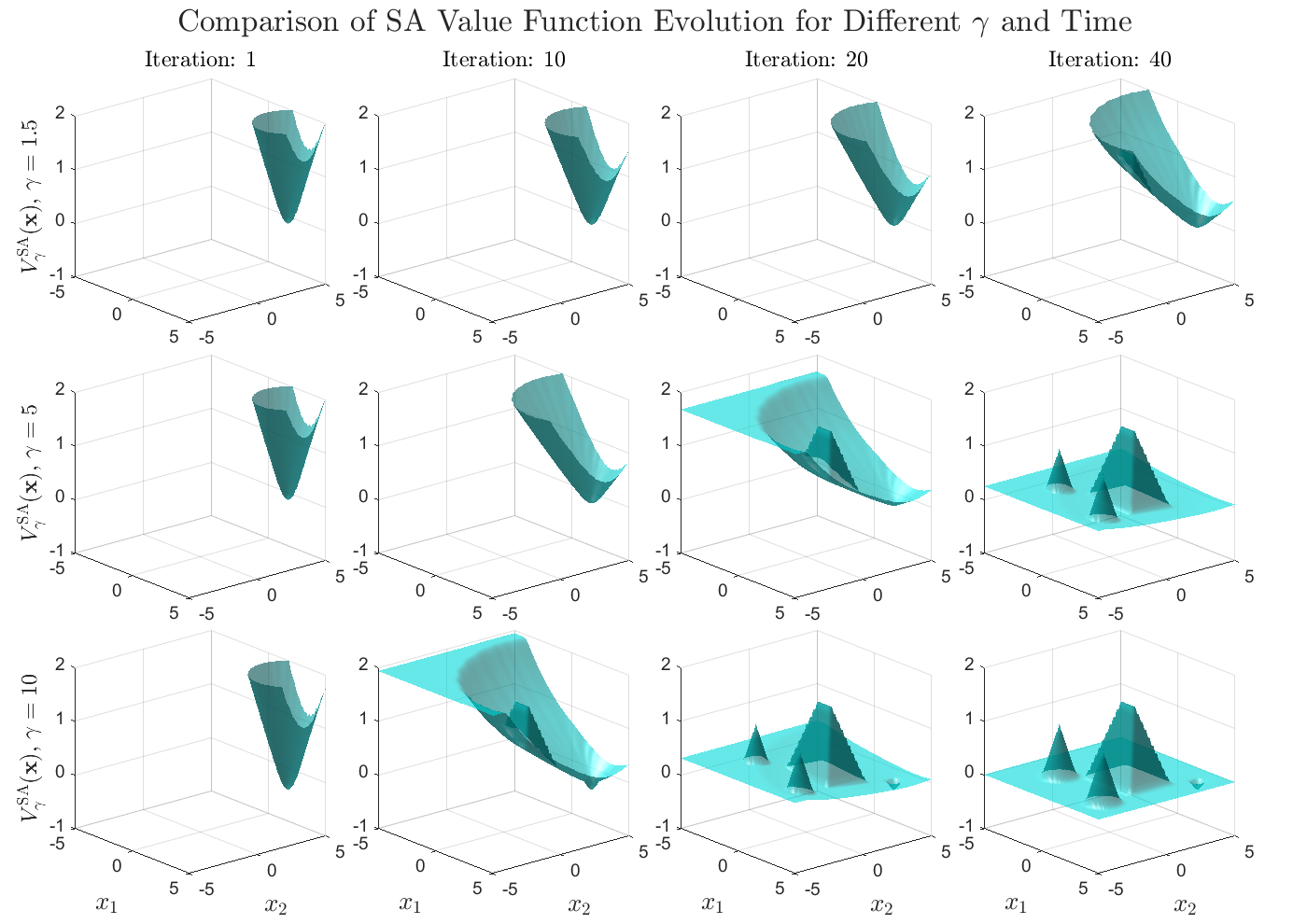}
    \caption{Evolution of the SA value functions over time for different $\gamma$. Top row, $\gamma = 1.5$. The SA value function evolves very slowly, and after 40 iterations, it is still far away from convergence. Medium row, $\gamma = 5$. The SA value function converges within 40 iterations. Bottom row, $\gamma = 10$. The SA value function converges within 20 iterations.} % Here, the values in the parentheses are the real-world computation time.}
 %   \vspace{-1em}
    \label{fig:SA_V_time_variations.png}
   \vspace{-1em}
\end{figure*}

\begin{figure}%[t]
    \centering
    \vspace{.5em}
    \includegraphics[width=0.49\textwidth]{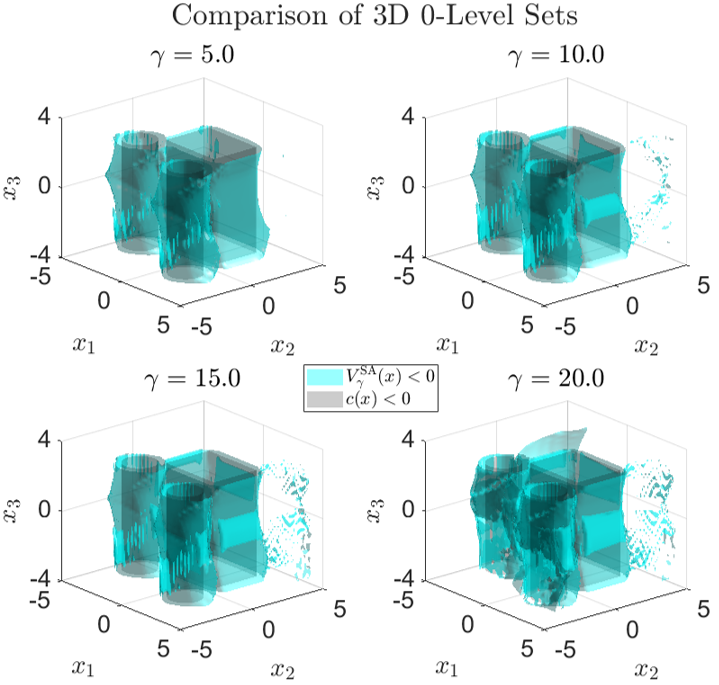}
    \caption{Zero level sets of constraint function $c(x)$ and SA value function $\VSA(x)$ (in the original 3D state space) for large $\gamma=5, 10, 15, 20$. As $\gamma$ increases, the characterization of zero level sets of $\VSA(x)$ becomes inaccurate.}
    \vspace{-1em}
    \label{fig:3D_level_Sets.png}
  % \vspace{-1em}
\end{figure}

\begin{figure}%[t]
    \centering
   %\vspace{.5em}
    \includegraphics[width=0.49\textwidth]{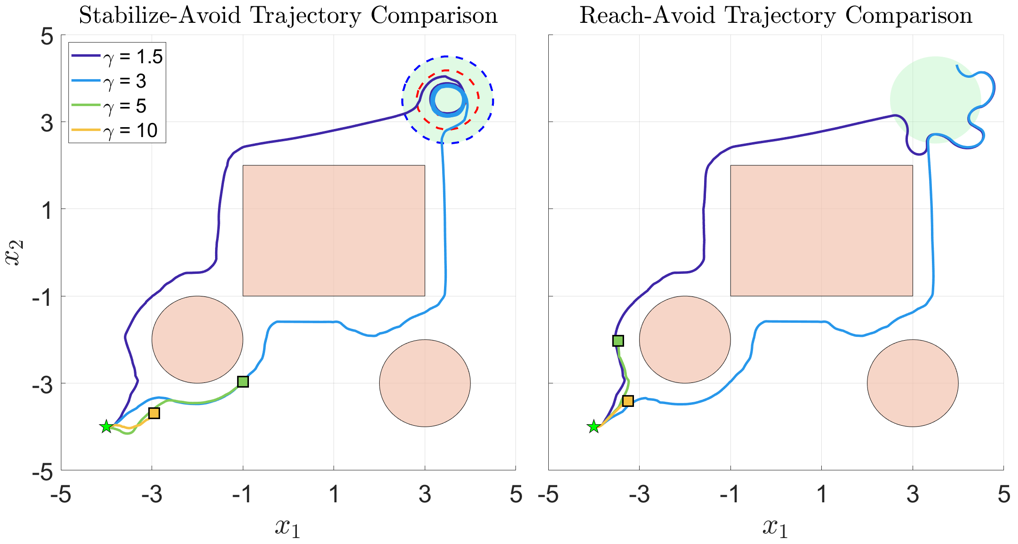}
    \caption{SA (left) and RA (right) trajectories of the 3D Dubins car that starts at $[-4;4;0]$ and aims to reach the green disk (target set $\mathcal{T}$). Left: the car safely reaches $\mathcal{I_\text{M}}$ (blue dashed lines) using $\pi_{RA}$ from~\eqref{eqn:RA_ctrl} and then stabilizes to $\mathcal{I_{\text{m}}}$ (the torus region enclosed by red dashed lines) using $\pi_{H}$ from~\eqref{eqn:SA_ctrl}. Right: the car safely reaches $\mathcal{T}$ using $\pi_{RA}$ from~\eqref{eqn:RA_ctrl}, but does not remain in the target set. %Since $\pi_{RA}$ is not defined for states inside the target, once we reach the target, we use the same control signal as it first reach the target. With this constant control signal, the trajectory leaves the target, and will be forced back after leaving. 
    Instead, the trajectory will repeatedly leave and return. The difference between trajectories is mainly caused by $\pi_{RA}$, which is derived by~\eqref{eqn:RA_ctrl}. However, there is no evidence that a larger $\gamma$ results in reaching the target faster. Note that only in the cases of $\gamma=1.5, 3$, the SA and RA tasks are accomplished. For larger $\gamma = 5, 10$, the trajectories terminate (square markers) before they can reach the target, due to insufficient snapshots of intermediate value functions $W(x,t)$ to synthesize the optimal control. }
    % The trajectories using $\pi_{RA}$ and $\pi_{H}$ are in magneta and cyan, respectively Left: the car safely reaches $\mathcal{I_\text{M}}$ (blue dashed lines) using RA controller~\eqref{eqn:RA_ctrl} and then stabilize to $\mathcal{I_{\text{m}}}$ (the torus region enclosed by red dashed lines). Right: the car safely reaches $\mathcal{T}$, but reapetedly leaves it after reaching.}
    \vspace{-1em}
    \label{fig:traj_3d_g=comb+0.2.png}
    %\vspace{-1em}
\end{figure}

In this section, we use a 3D Dubins car example to validate our theory. We empirically investigate the claims made in Section \ref{subsec:impact of gamma} regarding the impact of the discount factor, $\gamma$. We demonstrate that while any $\gamma > L_f$ is theoretically valid for recovering the reach-avoid set, its magnitude presents a critical trade-off between computational efficiency and the numerical fidelity required for controller synthesis. We will first show the benefits of a large $\gamma$ (Table.~\ref{tab:computation_times} and Fig.~\ref{fig:SA_V_time_variations.png}),  and then reveal its drawback (Fig.~\ref{fig:3D_level_Sets.png} and Fig.~\ref{fig:traj_3d_g=comb+0.2.png}).

The system is given by 
\begin{equation}\label{eq: 3D Dubins}
    \dot{x_1} = v\cos{x_3} + d_1 \quad \dot{x_2} = v\sin{x_3} + d_2 \quad \dot{x_3} = u,
\end{equation}
where $v=1$, $u\in [-\pi, \pi]$, and $d_1, d_2\in [-0.2, 0.2]$. The target set $\target=\{x|\ell(x)<0\}$, where $\ell (x) = (x_1-3.5)^2+(x_2-3.5)^2-1$, and the constraint set $\constraint = \{ x| \min(c_1(x),c_2(x),c_3(x)) < 0\}$, where $c_1(x) = 1 - (x_1+2)^2 - (x_2+2)^2$, $c_2(x) = 1 - (x_1-3)^2 - (x_2+3)^2 $, and $c_3(x) = 1 - \max\bigl( \tfrac{|x_1-1|}{2}, \, \tfrac{|x_2-\frac{1}{2}|}{1.5} \bigr)$. The Lipschitz constant for this system is $ L_f = 1$.

All simulations are conducted in MATLAB. The grid size used for computing $V_\gamma$ is $71\times 71\times 31$, and we summarize the computation time (averaged over $20$ runs) in Table \ref{tab:computation_times} for various $\gamma > L_f = 1$. It can be seen that as $\gamma$ increases from $1.5$ to $10^4$, the computation time decreases from $75.83s$ to $1.9s$. Fig. \ref{fig:SA_V_time_variations.png} provides a visual confirmation of this phenomenon, showing the evolution of the SA value function over time. For $\gamma=1.5$ (top row), the function evolves very slowly. In contrast, for $\gamma=10$ (bottom row), the function converges to its steady-state shape within just 20 iterations. These simulation results empirically validate our theoretical finding that a larger $\gamma$ strengthens the contraction of the Bellman backup operator, thereby speeding up convergence.

\begin{table}[b!]
\centering
% --- 组合拳：同时缩小字体和列间距 ---
\small % 1. 将表格内所有字体缩小
\setlength{\tabcolsep}{3pt} % 2. 极限压缩列间距 (可尝试2.5pt到4pt)
% -----------------------------------------
\caption{Average Computation Times of the RA value function for Different $\gamma$ Values}
\label{tab:computation_times}
\begin{tabular}{cr|cr|cr|cr}
\toprule
$\gamma$ & Time (s) & $\gamma$ & Time (s) & $\gamma$ & Time (s) & $\gamma$ & Time (s) \\
\midrule
1.5   & 75.83 & 10    & 17.96 & 100   & 4.47 & 1000  & 2.80 \\
3     & 43.93 & 20    & 10.56 & 200   & 3.55 & 5000  & 2.30 \\
5     & 29.55 & 50    & 6.01  & 500   & 3.02 & 10000 & 1.90 \\
\bottomrule
\end{tabular}
\end{table}

\begin{figure}%[t]
    \centering
    \vspace{.5em}
    \includegraphics[width=0.49\textwidth]{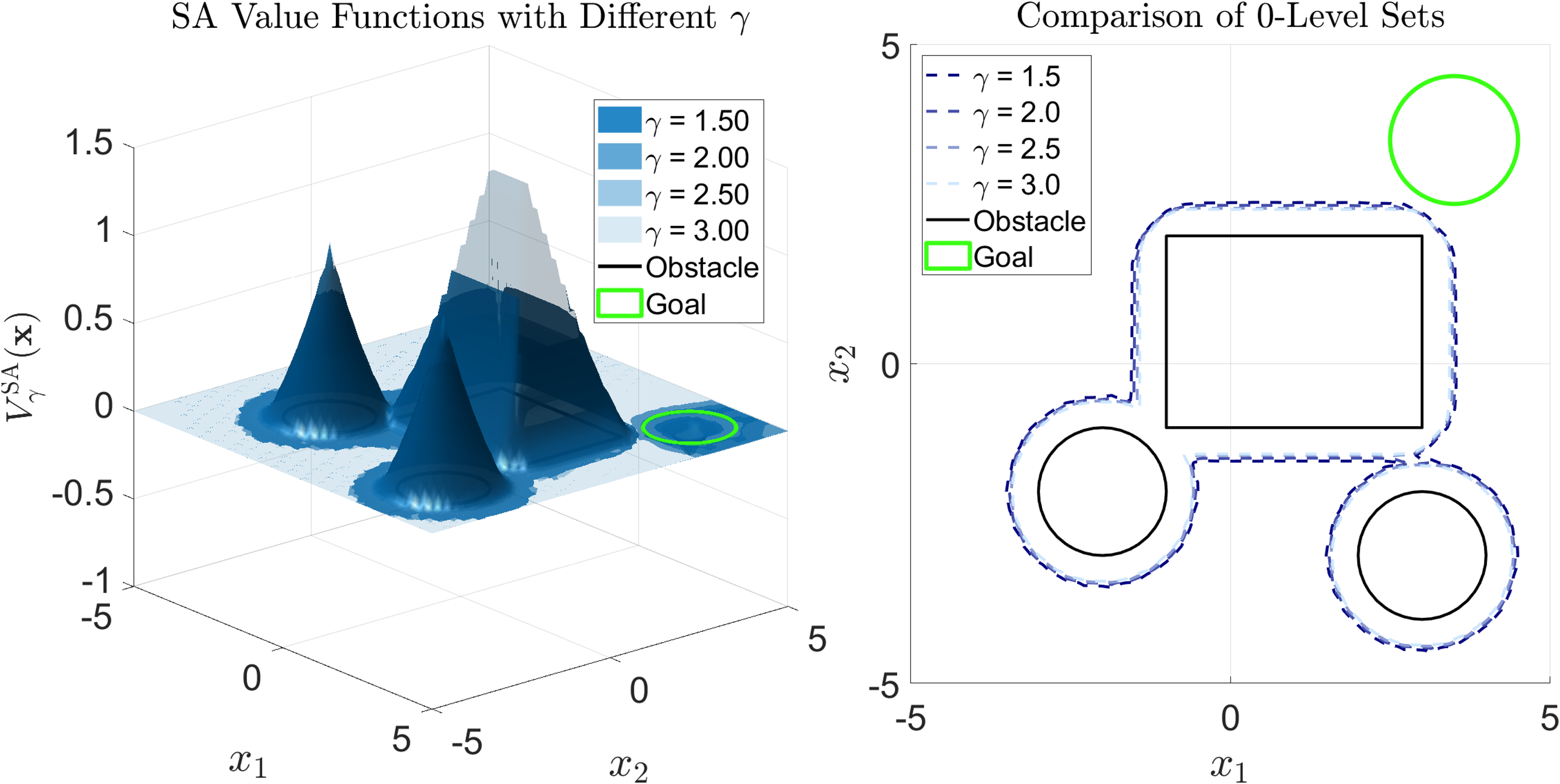}
    \caption{The projected SA value function $\VSA(x)$ and its zero level sets for moderate choices of $\gamma = 1.5, 2.0, 2.5, 3.0$. The regions enclosed by the green and black solid lines are the target $\target$ and obstacles, respectively, so the region outside the black solid lines is the constraint set $\constraint$. The regions enclosed by the dashed lines indicate the zero superlevel set of $\VSA(x)$, so the regions outside compose the $\mathcal{S} \mathcal{A}(\mathcal{T}, \mathcal{C})$ set. The SA value functions with different $\gamma$ provide identical zero level sets.}
    \vspace{-1em}
    \label{fig:V+Sets_g=1.5-3+0.2.png}
   %\vspace{-1em}
\end{figure}

Now, we verify that for appropriate choices of $\gamma$, the framework correctly identifies a consistent SA set. Fig. \ref{fig:V+Sets_g=1.5-3+0.2.png} displays the computed SA value functions, $V_{\gamma}^{\mathcal{SA}}(x)$, and their corresponding zero sublevel sets for $\gamma \in \{1.5, 2.0, 2.5, 3.0\}$. For visualization, we projected those SA functions to the $x_1-x_2$ plane. As predicted by Remark~1, while the shape of the value function changes with $\gamma$, the resulting certified SA set remains remarkably consistent, confirming that the underlying theoretical guarantee is robust to the choice of $\gamma$ within a reasonable range.

However, selecting an excessively large $\gamma$ introduces significant numerical artifacts. Fig. \ref{fig:3D_level_Sets.png} illustrates that as $\gamma$ increases to values like 5, 10, 15, and 20, the computed value function $V_{\gamma}^{SA}(x)$ begins to exhibit sharp oscillations near the obstacle boundaries. This leads to an inaccurate characterization of the zero level set, corroborating our analysis in Section \ref{subsec:impact of gamma}, i.e., large $\gamma$ values can cause numerical errors where the values of nearby states outside the SA set improperly influence the values in the SA set.

Another practical consequence of the numerical instability is the ``false convergence'' issue, which is most evident during control synthesis, as shown in Fig. \ref{fig:traj_3d_g=comb+0.2.png}. For well-chosen, smaller values of $\gamma$ (i.e., 1.5 and 3.0), the synthesized controllers successfully navigate the vehicle to the target region for both RA and SA tasks. In contrast, for larger values (i.e., $\gamma=5$ and $\gamma=10$), the generated trajectories terminate prematurely. This failure occurs because the rapid ``false convergence'' induced by a large $\gamma$ (as explained in Section \ref{subsec:impact of gamma}). The ``false convergence'' results in an insufficient number of time snapshots of the intermediate value function, $W(x,t)$, which are necessary to synthesize a complete, optimal control trajectory via~\eqref{eqn:RA_ctrl}.

In summary, the choice of $\gamma$ involves a clear and demonstrable trade-off. While a larger $\gamma$ offers substantial computational advantages, it can lead to numerical inaccuracies and premature termination of control synthesis. This underscores the importance of the practical guidelines suggested in Remark \ref{remark: effect of gamma} for selecting a value that balances computational tractability with the numerical precision required for generating valid controllers.

\section{CONCLUSIONS}
In this article, we presented a HJ-based framework to
solve both the infinite-horizon RA and SA games for general nonlinear continuous-time systems. We constructed a new discounted RA value function of which the zero sublevel set exactly characterizes the RA set. Moreover, we showed that the introduction of the discount factor leads to many desirable properties of the value function: the Lipschitz continuity under certain assumptions, the uniqueness of the viscosity solution to a corresponding HJI-VI, and the contraction property of the associated Bellman backup operator that guarantees convergence from arbitrary initializations. By integrating our RA strategy with R-CLVF, we developed a two-step framework to construct a SA value function, of which the zero sublevel set fully recovers the desired SA set. Finally, we provided controller synthesis approaches for both the RA and SA tasks.

\section*{ACKNOWLEDGMENT}
We thank Professor Donggun Lee from North Carolina State University, Jingqi Li
and Jason Choi from UC Berkeley, Haimin Hu from Princeton University, and Sander Tonkens, Will Sharpless, and Dylan Hirsch from UC San Diego for their insightful comments and valuable discussions.
%\vspace{1em}
\appendix
\subsection{Proof of the Uniqueness of the Viscosity Solution}
\begin{proof}
    Let $V^1, V^2$ be the sub- and supersolutions of \eqref{eq: HJI-VI}, respectively. % and assume by contradiction that there exists $\delta > 0$ and $\Tilde{x}$ such that $V^1(x) - V^2(x) = \delta$. 
    Let $||x||$ be the Euclidean norm of $x\in \R^n$.
    Consider \( \Phi : \mathbb{R}^{2N} \to \mathbb{R} \),
    \begin{align}\label{eq: def_Phi}
    \Phi(x,y) := V^1(x) - V^2(y) - \frac{||x - y||^2}{2\varepsilon}
    \end{align}
    where \( \varepsilon \) is a positive parameter to be chosen conveniently. Let us assume by contradiction that there is \( \delta > 0 \) and \( \Tilde{x} \) such that \( V^1(\Tilde{x}) - V^2(\Tilde{x}) = \delta \). Then we have, 
    \begin{align}\label{eq: uni_ctd}
    \frac{\delta}{2} < \delta = \Phi(\Tilde{x}, \Tilde{x}) \leq \sup \Phi(x,y).
    \end{align}
    Since \( \Phi \) is continuous and \(\lim_{|x| + |y| \to \infty} \Phi(x,y) = -\infty\) when \(x\neq y\), 
    there exist \( \bar{x}, \bar{y} \) such that
    \begin{align}\label{eq: def_xybars}
    \Phi(\bar{x}, \bar{y}) = \sup \Phi(x,y).
    \end{align}
    Thus, the inequality $\Phi(\bar{x}, \bar{x}) + \Phi(\bar{y}, \bar{y}) \leq 2\Phi(\bar{x}, \bar{y})$ holds, so we easily get
    \begin{align}\label{eq: (x-y)2}
    \frac{||\bar{x} - \bar{y}||^2}{\varepsilon} \leq V^1(\bar{x}) - V^1(\bar{y}) + V^2(\bar{x}) - V^2(\bar{y}).
    \end{align}
    Then the boundedness of \( V^1 \) and \( V^2 \) implies
    \begin{align}\label{eq: bound_bar_x-y}
    ||\bar{x} - \bar{y}|| \leq k\sqrt{\varepsilon}.
    \end{align}
    for a suitable constant \( k \) irrelevant of $x,y$. By plugging \eqref{eq: bound_bar_x-y} into \eqref{eq: (x-y)2} and using the uniform continuity of \( V^1 \) and \( V^2 \), we get
    \begin{align}\label{eq: bar_x-y_mod}
    \frac{||\bar{x} - \bar{y}||^2}{\varepsilon} \leq \omega(\sqrt{\varepsilon}),
    \end{align}
    for some modulus \( \omega \), i.,e., a function $\omega : [0, +\infty) \to [0, +\infty)$ that is continuous, nondecreasing, and satisfies \(\omega(0) = 0\). Next, define the \( C^1 \) test functions
    \begin{align*}
    \varphi(x) &:= V^2(\bar{y}) + \frac{||x - \bar{y}||^2}{2\varepsilon}, \\
    \phi(y) &:= V^1(\bar{x}) - \frac{||\bar{x} - y||^2}{2\varepsilon},
    \end{align*}  
    and observe that, by definition of \( \bar{x}, \bar{y} \), \( V^1 - \varphi \) attains its maximum at \( \bar{x} \) and \( V^2 - \phi \) attains its minimum at \( \bar{y} \). It is easy to compute 
    \begin{align}\label{eq: equal_D}
    D_x\varphi(\bar{x}) &= \frac{\bar{x} - \bar{y}}{\varepsilon} = D_y\phi(\bar{y})
    \end{align}  
    By definition of viscosity sub- and supersolution, we have
    \begin{align}
        \max\bigl\{&\min \{ H\left( \bar{x}, D_{x}\varphi(\bar{x}) \right) - \gamma V^1(\bar{x}), \notag\\
        & \ell(\bar{x})-V^1(\bar{x})\}, c(\bar{x})-V^1(\bar{x})\bigr\}\geq 0, \label{eq: uni_sub_ineq}\\
        \max\bigl\{&\min \{ H\left( \bar{y}, D_{y}\phi(\bar{y}) \right) - \gamma V^2(\bar{y}), \notag\\
        & \ell(\bar{y})-V^2(\bar{y})\}, c(\bar{y})-V^2(\bar{y})\bigr\}\leq 0, \label{eq: uni_sup_ineq}    
    \end{align}
    From \eqref{eq: uni_sup_ineq} we have
    \begin{align}\label{eq: uni_sup_c-V2}
        c(\bar{y})-V^2(\bar{y}) \leq 0,
    \end{align}
    and one of the following holds:
    \begin{align}
        H\left( \bar{y}, D_{y}\phi(\bar{y}) \right) - \gamma V^2(\bar{y}) \leq 0, \label{eq: uni_sup_H-V2}\\
        \ell(\bar{y})-V^2(\bar{y})\leq 0, \label{eq: uni_sup_l-V2}
    \end{align}
    From \eqref{eq: uni_sub_ineq} we have
    \begin{align}\label{eq: uni_sub_c-V1}
        c(\bar{x})-V^1(\bar{x}) \geq 0,
    \end{align}
    or both of the following hold:
    \begin{align}
        H\left( \bar{x}, D_{x}\varphi(\bar{x}) \right) - \gamma V^1(\bar{x}) \geq 0, \label{eq: uni_sub_H-V1}\\
        \ell(\bar{x})-V^1(\bar{x})\geq 0, \label{eq: uni_sub_l-V1}
    \end{align}
    Let us first assume that \eqref{eq: uni_sub_c-V1} and \eqref{eq: uni_sup_c-V2} are true. After rearranging the terms and using \eqref{eq: bound_bar_x-y}, we get
    \begin{align}\label{eq: V1-V2_c}
        V^1(\bar{x}) - V^2(\bar{y}) &\leq  c(\bar{x}) - c(\bar{y}) \leq L_c||\bar{x} -\bar{y}|| \leq L_c k\sqrt{\varepsilon},
    \end{align}
    which implies that $\Phi(\bar{x}, \bar{y})\leq L_c k\sqrt{\varepsilon}$. Similarly, if \eqref{eq: uni_sub_l-V1} and \eqref{eq: uni_sup_l-V2} hold, we can show that
    \begin{align}\label{eq: V1-V2_l}
        V^1(\bar{x}) - V^2(\bar{y}) &\leq  \ell(\bar{x}) - \ell(\bar{y}) \leq L_\ell||\bar{x} -\bar{y}|| \leq L_\ell k\sqrt{\varepsilon}.
    \end{align}

    Finally, assuming \eqref{eq: uni_sup_H-V2} and \eqref{eq: uni_sub_H-V1} hold, we get
    \begin{align}\label{eq: V1-V2_H_1}
        V^1(\bar{x}) - V^2(\bar{y}) &\leq \frac{1}{\gamma} \bigl(H\left( \bar{x}, D_{x}\varphi(\bar{x}) \right) - H\left( \bar{y}, D_{y}\phi(\bar{y}) \right)\bigr).
    \end{align}
    By the compactness of $\cset$ and $\dset$ and the Lipschitz continuity of $f$ in $x$, the Hamiltonian $H$ has the property that for a fixed $p\in \R$
    \begin{equation}\label{eq: H_Lip}
        |H(x,p) - H(y,p)|\leq |p|L_f||x - y||,
    \end{equation}
    for any $x,y\in \R$. Plugging \eqref{eq: H_Lip} back to \eqref{eq: V1-V2_H_1} and invoking \eqref{eq: equal_D} and \eqref{eq: bar_x-y_mod}, we get
    \begin{align}\label{eq: V1-V2_H_2}
        V^1(\bar{x}) - V^2(\bar{y}) &\leq \frac{L_f}{\gamma}\omega(\sqrt{\varepsilon}).
    \end{align}
    Combining \eqref{eq: V1-V2_c}, \eqref{eq: V1-V2_l}, \eqref{eq: V1-V2_H_2} and \eqref{eq: def_Phi}, we obtain
    \begin{align}
        \Phi(\bar{x}, \bar{y})\leq \frac{1}{3}(L_c k \sqrt{\varepsilon} +L_\ell k \sqrt{\varepsilon} + \frac{L_f}{\gamma}\omega(\sqrt{\varepsilon})),
    \end{align}
    and the right-hand side can be made smaller than $\frac{\delta}{2}$ for $\varepsilon$ small enough, a contradiction to \eqref{eq: uni_ctd} and \eqref{eq: def_xybars}. To conclude, we have proven that a subsolution of \eqref{eq: HJI-VI} will not be larger than a supersolution. Since $V_{\gamma}$ is both a sub- and supersolution, it is the unique one.
\end{proof}

\subsection{Derivation of Equation~\texorpdfstring{\eqref{eq:ex_RA_val}}{(44)}}
The solution of~\eqref{eq:ex_sys} under control~\eqref{eq:ex_ctrl} is:
\begin{align} \label{eq:ex_sol}
    \xi_x^u (t) = \begin{cases}
        0 & |x| \leq t \\ x - \text{sgn}(x)\cdot t& 0\leq t \leq |x|
    \end{cases}.
\end{align}
Let us first find $V_\gamma(0)$. Since $\xi_{x=0}^u (t) = 0$, plugging in~\eqref{eq: def-inf}:
\begin{align*}
    V_\gamma(0) &= \inf_{t\in[0,\infty)} \max \bigl \{ e^{-\gamma t} (|0|-1), \max_{s \in [0,t] }e^{-\gamma s} (|0|-2) \bigr \},\\
    & = \inf_{t\in[0,\infty)} \max \bigl \{- e^{-\gamma t} , -2e^{-\gamma t}  \bigr \}, \\
    &= \inf_{t\in[0,\infty)} (- e^{-\gamma t}) = -1.
\end{align*}
We next derive $V(x)$ for $x>0$. Use~\eqref{eq: def-inf}: \begin{align*}
    &V_\gamma(x) \\ 
    =&  \inf_{t\in[0,\infty)} \max \bigl \{ e^{-\gamma t} (|\xi_x^u (t)|-1), \max_{s \in [0,t] }e^{-\gamma s} (|\xi_x^u (s)|-2) \bigr \}.
\end{align*}
Denote $A(s;\gamma) := e^{-\gamma s} (|\xi_x^u (s)|-2) $, we have: \begin{align*}
    A(s;\gamma) = \begin{cases}
        e^{-\gamma s} (|x-s|-2) & s\leq x \\
        -2 e^{-\gamma s} & s>x
    \end{cases}.
\end{align*} Take the derivative: \begin{align*}
    \dot A(s;\gamma) = \begin{cases}
        \gamma e^{-\gamma s} (2+s-x-\frac{1}{\gamma} ) & s\leq x \\
        2\gamma e^{-2\gamma s}  & s>x
    \end{cases}.
\end{align*} When $s>x$, it can be seen that $\dot A(s;\gamma)  > 0$ for any $\gamma >0$. When $s\leq x$, $\dot A(s;\gamma)$ can be negative for some $\gamma$: eg., $\gamma >\frac{1}{2}$ and $s<x + \frac{1}{\gamma}-2 < x$. In conclusion, for any $\gamma$, the sign of $\dot A(s;\gamma)$ either doesn't change or changes only one time from negative to positive. Therefore, $\max_{s \in [0,t] } A(s;\gamma)$ attains the maximum value either at $s=0$ or $s = t$:\begin{align*}
    \max_{s \in [0,t] } A(s;\gamma) = \max ( x-2 , e^{-\gamma t} (|\xi_x^u (t)|-2)  ).
\end{align*} The value function becomes \begin{align*}
    &V_\gamma(x) \\ 
    =&  \inf_{t\in[0,\infty)} \max \bigl \{ e^{-\gamma t} (|\xi_x^u (t)|-1), x-2 , e^{-\gamma t} (|\xi_x^u (t)|-2)  \bigr \}, \\
    =&\inf_{t\in[0,\infty)} \max \bigl \{ e^{-\gamma t} (|\xi_x^u (t)|-1), x-2   \bigr \} \\
    =& \max \bigl \{ \inf_{t\in[0,\infty)} e^{-\gamma t} (|\xi_x^u (t)|-1), x-2   \bigr \}.
\end{align*} Notice that up to this step, $\gamma$ have no impact on the value function. This is because no matter what is the value of $\gamma$, the this inequality always hold: $$e^{-\gamma t} (|\xi_x^u (t)|-2)< e^{-\gamma t} (|\xi_x^u (t)|-1).$$ In the next step, we will see how $\gamma$ changes the value function.

Now, inspect $\inf_{t\in[0,\infty)} e^{-\gamma t} (|\xi_x^u (t)|-1)$. We have: \begin{align*}
    &\inf_{t\in[0,\infty)} e^{-\gamma t} (|\xi_x^u (t)|-1) \\
    &= \min ( \min_{t\in[0,x] }  e^{-\gamma t} (|x-t|-1), \inf_{t\in[x,\infty)} - e^{-\gamma t} ) \\
    & = \min ( \min_{t\in[0,x] }  e^{-\gamma t} (|x-t|-1), -e^{-\gamma x} )
\end{align*} denote $B(t;\gamma) := e^{-\gamma t} (|x-t| -1 )$, $t\in[0,x]$. We have: \begin{align*}
    B(t;\gamma) = 
        e^{-\gamma t} (x-t-1).
\end{align*} Take the derivative: \begin{align*}
    \dot B(t;\gamma) = 
        \gamma e^{-\gamma t} (1+t-x-\frac{1}{\gamma} ) 
\end{align*} For $\gamma \leq 1$, $\dot B(t;\gamma)<0$ for all $t\in[0,x]$. The minimum is attained at $t = x$: \begin{align*}
    \min_{t\in[0,x]}B(t;\gamma) = 
        -e^{-\gamma x}.
\end{align*} Therefore, we have \begin{align*}
     V_\gamma(x)
    =& \max \bigl \{ \inf_{t\in[0,\infty)} e^{-\gamma t} (|\xi_x^u (t)|-1), x-2   \bigr \} \\
    =& \max \bigl \{ \min ( -e^{-\gamma x}, -e^{-\gamma x} ), x-2   \bigr \} \\
    =& \max \bigl \{  -e^{-\gamma x} , x-2   \bigr \}
\end{align*} 
For $\gamma > 1$, when $x<1-\frac{1}{\gamma}$, $\dot B(t;\gamma)$ is positive for all $t\in[0,x]$, therefore the minimum is attained at $t=0$. We have \begin{align*}
     V_\gamma(x)
    = \max \bigl \{ \min ( x-1, -e^{-\gamma x} ), x-2   \bigr \}. 
\end{align*} When $x<1-\frac{1}{\gamma}$ $\dot B(t;\gamma)$ can be negative when $0\leq t<x + \frac{1}{\gamma}-1 < x$. But the sign of $\dot B(t;\gamma)$ changes only once from negative to positive at $t = x + \frac{1}{\gamma} - 1$. Therefore, the infimum of $B(t;\gamma)$ is attained at $t = x + \frac{1}{\gamma} - 1$: \begin{align*}
    &\inf_{t\in[0,x]} e^{-\gamma t} (|\xi_x^u (t)|-1) \\
    = &  e^{-\gamma (x+\frac{1}{\gamma} -1)} (| \xi_x^u (x+\frac{1}{\gamma} -1) |-1 )\\
    =& e^{-\gamma x -1 + \gamma} (1-\frac{1}{\gamma}-1) = -\frac{1}{\gamma} e^{-\gamma x -1 + \gamma}
\end{align*} Therefore, we have \begin{align*}
     V_\gamma(x)
    =& \max \bigl \{ \inf_{t\in[0,\infty)} e^{-\gamma t} (|\xi_x^u (t)|-1), x-2   \bigr \} \\
    =& \max \bigl \{ \min ( -\frac{1}{\gamma} e^{-\gamma x -1 + \gamma}, -e^{-\gamma x} ), x-2   \bigr \} 
\end{align*} 
Since the solution~\eqref{eq:ex_sol} is symmetric, the value for $x<0$ can be obtained by adding a negative sign before $x$. Combined, we get~\eqref{eq:ex_RA_val}.
%%%%%%%%%%%%%%%%%%%%%%%%%%%%%%%%%%%%%%%%%%%%%%%%%%%%%%%%%%%%%%%%%%%%%%%%%%%%%%%%%%%%%%%%%%%%%%%%%%%%%%%%%%%%%%%%%%%%%%%%%%%%%%%%%%%%%%%%%%%%%%%%%%%%%%%%%%%%  
\bibliographystyle{IEEEtran}
\bibliography{reference}
\vspace{-2em}
\begin{IEEEbiography}
[{\includegraphics[width=1in,height=1.25in,clip,keepaspectratio]{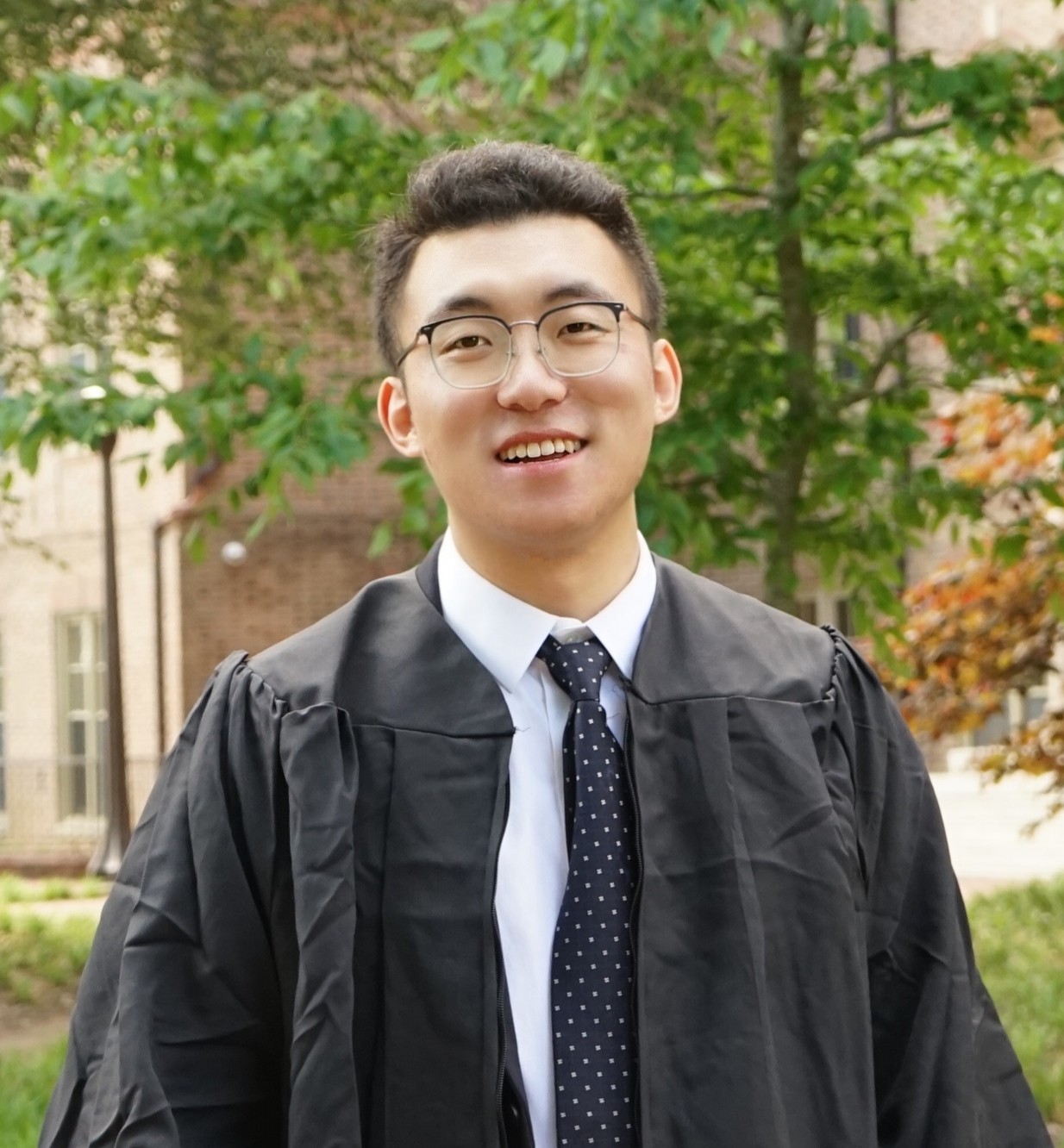}}]{Boyang Li}{\,}(bol025@ucsd.edu) received a B.S. degree in Mathematics and Physics (double major) from William \& Mary in May 2022. He is currently pursuing a Ph.D. degree in Mechanical Engineering at UC San Diego. Boyang is interested in developing control methodologies that possess both mathematical formality and versatility in application to various robotics systems, drawing ideas from optimal control, nonlinear systems, and deep learning.
\end{IEEEbiography}
\vspace{-2em}
\begin{IEEEbiography}
[{\includegraphics[width=1in,height=1.25in,clip,keepaspectratio]{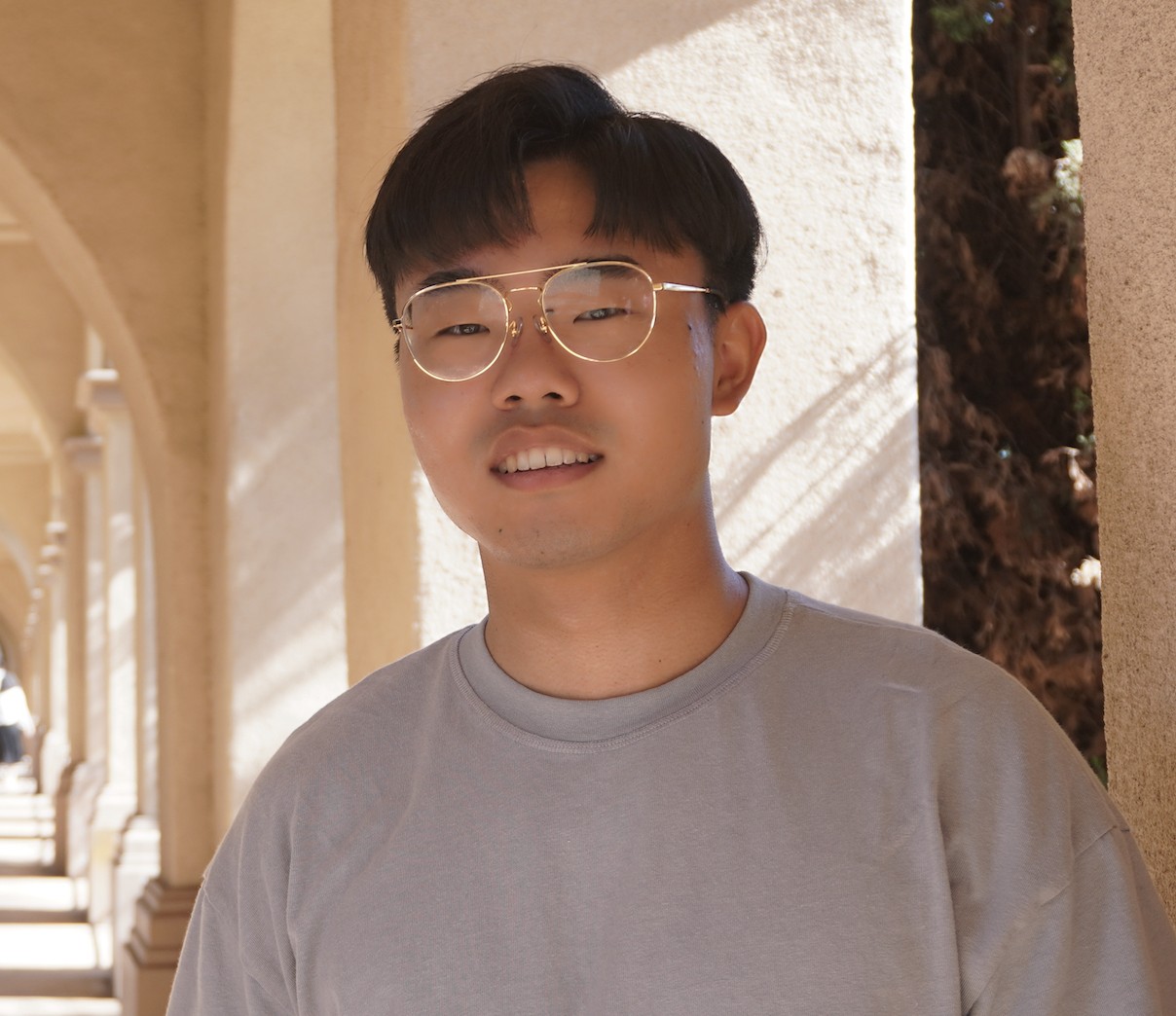}}]{Zheng Gong}{\,}(zhgong@ucsd.edu) received a B.S. degree in Vehicle Engineering from Beijing Institute of Technology in June 2019, and a M.S. degree in Mechanical Engineering and Materials Science from Duke University in May 2021. He is currently pursuing a Ph.D. degree in Mechanical Engineering at UC San Diego. Zheng is interested in finding connections between CLF and HJR analysis, developing efficient computation methods, and exploring data-driven approaches to solving PDEs.
\end{IEEEbiography}
%\vspace{-45em}
\begin{IEEEbiography}
[{\includegraphics[width=1in,height=1.25in,clip,keepaspectratio]{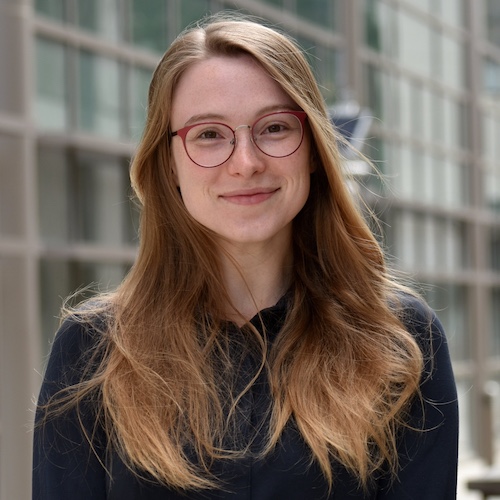}}]{Sylvia L. Herbert}{\,}(sherbert@ucsd.edu) is an Assistant Professor in Mechanical and Aerospace Engineering at the University of California San Diego.  Prior to joining UCSD, she received her PhD in Electrical Engineering from UC Berkeley, where she studied with Professor Claire Tomlin on safe and efficient control of autonomous systems.  Before that she earned her BS/MS at Drexel University in Mechanical Engineering.  She is the recipient of the ONR Young Investigator Award, Hellman Fellowship, UCSD JSOE Early Career Faculty Award, UC Berkeley Chancellor’s Fellowship, NSF GRFP, UC Berkeley Outstanding Graduate Student Instructor Award, and the Berkeley EECS Demetri Angelakos Memorial Achievement Award for Altruism. 
\end{IEEEbiography}
\end{document}